# CIR bridge for modeling of fish migration on sub-hourly scale

Hidekazu Yoshioka[a,*]

[a] Japan Advanced Institute of Science and Technology, 1-1 Asahidai, Nomi 923-1292, Japan

[*] Corresponding author: yoshih@jaist.ac.jp, ORCID: 0000-0002-5293-3246

**Abstract**

Bridges, which are stochastic processes with pinned initial and terminal conditions, have recently been applied to various problems. We show that a bridge based on the Cox–Ingersoll–Ross process, called a CIR bridge in this paper, reasonably models the intraday number of migrating fish at an observation point in a river. The studied fish migrates between sunrise and sunset each day, which are considered the initial and terminal times, respectively. The CIR bridge is well-defined as a unique pathwise continuous solution to a stochastic differential equation with unbounded drift and diffusion coefficients and potentially represents the on–off intermittency of the fish count data. Our bridge is theoretically novel in that it admits closed-form time-dependent averages and variances, with which the model parameters can be identified efficiently, and is computable by a recently-developed one-step numerical method. The CIR bridge is applied to the sub-hourly migration data of the diadromous fish *Plecoglossus altivelis altivelis* in the Nagara River, Japan, from February to June.

***Statements & Declarations***

**Acknowledgments:** The author would like to express his gratitude to Japan Water agency, Nagara River Estuary Barrage Office for providing the valuable fish migration data. The author thanks the two reviewers for their helpful comments and suggestions.

**Funding:** This study was supported by the Japan Science and Technology Agency (PRESTO No. JPMJPR24KE).

**Competing Interests:** The authors have no relevant financial or non-financial interests to disclose.

**Data Availability:** The data will be made available at reasonable request to the corresponding author.

**Declaration of Generative AI in Scientific Writing:** The authors did not use generative AI for the scientific writing of this manuscript.

**Contribution:** all parts of this manuscript were prepared by the author.

## 1. Introduction

### 1.1 Background

Fishes migrating between a river and a sea, called diadromous, are key drivers in the food web and nutrient cycling [1,2]. Diadromous fishes include high-demand species, such as salmonids and sturgeons [3] and lampreys [4]. Their populations have been declining globally due to the habitat fragmentation and degradation caused by anthropogenic activities such as dam construction, and their conservation through countermeasures such as restocking and fisheries is an urgent issue [5,6,7]. The population dynamics of diadromous fish species should therefore be evaluated for the planning of their conservation strategies.

The population dynamics of diadromous fishes can be quantified by counting them at some observation point in a river where they migrate. Their daily count data have been used to assess the effects of climate change and river hydrology on the migration of Atlantic salmon [8], Chinook salmon [9], and European eel [10,11]. The effects of anthropogenic interventions on fish migration were investigated, with a focus on the trapping of Atlantic salmon, using daily fish count data [12]. A passive counting technology with statistical uncertainty estimation was developed and applied to the migration of Sockeye salmon [13]. Daily environmental DNA sampling was recently implemented to complement sonar-based fish count data [14].

Although daily count data are valuable for long-term fish observation and management, fish migration phenomena involve sub-daily timescales. Such finer-resolution data are needed for a better understanding of the biology of fish migration and for the design of effective fish conservation strategies. The diel upstream migration of major diadromous fish species, including salmonids, steelhead, and lampreys, shows that the peak time of the migrating fish count varies across species [15]. The diurnal and tidal cycles affect the migration behavior of diadromous fish species, suggesting the critical influence of operating tidal gates on their population dynamics [16]. Bozeman et al. [17] reviewed the effects of sub-daily flow regulation for hydropower generation on the migration of diadromous fishes. Norman et al. [18] emphasized the importance of understanding the fine temporal dynamics of fish migration for the safe operation of water infrastructure. Hourly observations have shown that juvenile Ayu *Plecoglossus altivelis altivelis* (*P. altivelis*) in Japan migrates upstream along a river during daytime [19,20]. Despite the information that can be obtained from fine descriptions of fish migration for a better understanding of their biology and ecology and for their sustainable coexistence with humans, such observations remain scarce in literature. A sub-daily scale mathematical model for fish migration will form a part of the foundation of the solution to this issue.

Mathematical models for fish migration are often stochastic, possibly because of the difficulty of fully describing the whole mechanisms involved in fish migration. The density-dependent population dynamics of migratory fishes were discussed using a projection matrix model driven by early-life stochasticity [21]. An empirical Bayesian approach was incorporated into a Markov chain model of spatial fish movement across multiple habitats [22]. A multiple regression model was developed for assessing the survival and migration timing of reintroduced fishes [23]. Cherbero et al. [24] examined the use of linear mixed-effect models to evaluate the upstream migration timing of European shads for spawning. The effects of water temperature and hydrology on the migration strategies of the bull trout were analyzed using a linear spatial capture–recapture model that concurrently studies the migration and survival of radio-tagged fishes [25]. Stochastic individual-based models that regard individual fishes as diffusive particles were utilized to assess ecologically friendly reservoir operation [26]. Life-cycle models incorporating stochastic components representing the uncertainties involved in environmental and population dynamics have been adopted for the management of endangered species [27] and invasive species [28]. Torstenson and Shaw [29] proposed a conceptual model for explaining the seasonal dependence of the migration timing of migratory animals and subjected it to a sensitivity analysis. Yoshioka and Yamazaki [30] described the daily time series of the fish count data of *P. altivelis* in a river as a continuous-time stochastic process driven by self-exciting jumps. Yoshioka [31] applied the superposition of jump-driven affine processes to similar data of the same fish species in multiple rivers in Japan. Sato and Seguchi [20] studied 10 min migration data and analyzed its relation to water quantity and quality but did not conduct mathematical modeling. Mathematical models that deal with the fine (sub-daily) dynamics of migrating fish populations have not been studied despite their potential importance in both theory and application. We address this problem in this paper by obtaining fish migration data.

### 1.2 Aim and contribution

This paper aims to present a novel mathematical model for fish migration on a sub-hourly scale. We stochastically model the diurnal upstream migration behavior of juveniles of *P. altivelis*, a major inland fishery resource in Japan (e.g., Tsukamoto et al. [32]). According to the Fisheries Agency[1], this fish accounts for approximately half the value of inland fishery production in recent years in Japan. *P. altivelis* has a typical lifespan of one year, and we investigate its upstream migration from seas to a river in Japan from February to June, the juvenile stage. This fish migrates only during daytime each day, suggesting that the sunrise and sunset times can be considered the initial and terminal times of fish migration. Our state variable is the fish count data at a fixed observation point in a river.

A stochastic process with a pinned terminal value is called a bridge, and a bridge in the form of a diffusion-driven stochastic differential equation (SDE) is called a diffusion bridge. Diffusion bridges have been used in various research fields, such as areas involving asset valuation, e.g., finance and economics [33,34,35], physics and chemistry [36,37,38], animal movement paths [39,40,41], random graph generation [42], and machine learning [43,44,45]. Advanced bridges have also been investigated, such as Brownian bridges with resetting [46,47] and diffusion bridges with area constraints [48]. In our context, a time series of daily fish count data between sunrise and sunset is considered a bridge—a novel viewpoint about fish migration. Here, we should explicitly state that a "bridge" in this paper is a stochastic process with a prescribed terminal condition, which is a broader concept than those based on the conditioning reviewed in the next paragraph.

[1] https://www.jfa.maff.go.jp/j/enoki/attach/pdf/naisuimeninfo-41.pdf (Last accessed on May 14, 2025)

Diffusion bridges are often formulated through conditioning at the terminal time (often called h-transform), where an additional nonlinear drift term due to the conditioning arises in the corresponding SDEs [49,50,51]. An energy minimization approach has also been used to formulate (diffusion and jump) bridges, where the additional drift is interpreted as a control variable [52,53,54]. In the use of a bridge to model a nonnegative state variable, such as the fish count in a unit time, the underlying process should also be almost surely (a.s.) nonnegative. The Cox–Ingersoll–Ross (CIR) process is the simplest such SDE found in economics [55]; having efficiently computable statistics and a unique solution that is a.s. nonnegative, it has been applied to many problems, such as disease propagation [56], nurse deployment optimization [57], water quality dynamics [58], and heat engine analysis [59].

A bridge based on the CIR process is called a CIR bridge in this paper. CIR bridges based on the h-transform have been applied to option pricing [60], business scenario simulations of stocks [61], and hysteresis modeling for magnetic materials [62]. We apply a CIR bridge to the upstream migration of *P. altivelis*, where the state variable is the fish count in a unit time observed at a fixed point in a river. Our CIR bridge differs from classical ones in that it depends not on conditioning but on a more heuristic formulation to generate sample paths pinned at the terminal (and initial) time. With this formulation, we avoid assuming a sufficiently small diffusion constant, also called volatility—a common restriction regarding parameter ranges in the literature [63].

Another advantage of the proposed model is that its form enables modern numerical methods for simulating CIR-type processes, called integrated variance implicit (iVi) schemes [64,65], to be applied with slight modifications. Classical numerical methods for CIR and related diffusion processes with non-Lipschitz diffusion coefficients do not cover high-volatility cases (e.g., loss of nonnegativity of numerical solutions) [66,67,68]. However, bridges usually have drift coefficients that become unbounded near the terminal time, causing the assumptions in classical numerical methods to break down. The iVi scheme avoids this issue and is robust throughout the range of volatility values.

We apply the CIR bridge to the unique 10 min count data of *P. altivelis* in the Nagara River, Japan. The data suggest that migration occurs during daytime, as indicated by bursts in fish count. The tractable nature of the CIR bridge, particularly the closed-form availability of its time-dependent average and variance, is used for parameter fitting. Daily fish count data obtained at the same observation point are then used to investigate the seasonal dependence of the model parameters in the CIR bridge. The fitted model is investigated in terms of statistical moments and pathwise properties, with the burst behavior of sample paths being computationally investigated in particular. Our model only considers temporal fish migration dynamics although it only characterizes a part of ecosystem dynamics that depend both on space and time (e.g., Spagnolo et al. [69]). Nevertheless, the data analysis of the fish count and water temperature gives some hints toward advanced mathematical modeling of fish migration.

The modeling and analysis conducted in this paper show that the CIR bridge is a well-posed stochastic process model, and it can be applied to the upstream migration of *P. altivelis*, which is a phenomenon characterized by highly stochastic and intermittent behavior. This paper thus contributes to both the modeling and application of the new mathematical model for describing fish migration.

The rest of this paper is organized as follows. **Section 2** presents the CIR process and bridge. **Section 3** mathematically analyzes the bridge, focusing on its statistical well-posedness. **Section 4** presents the application of the proposed CIR bridge to the sub-hourly fish count data. **Section 5** concludes this paper and presents study prospects. **Appendix** contains proofs of the propositions (**Section A.1**) and auxiliary results (**Sections A.2 and A.3**).

## 2. Mathematical model

As usual in continuous-time stochastic modeling with SDEs (e.g., Chapter 1 in Capasso and Bakstein [70]), we consider a complete probability space $(\Omega, \mathbb{F}, \mathbb{P})$. Here, $\Omega$ is the sample space (containing all possible outcomes), $\mathbb{F}$ is the event space (containing all events), and $\mathbb{P}$ is the probability function. Time $t$ is a continuous parameter.

### 2.1 CIR process

We briefly review the classical CIR process following Alfonsi [71]. All descriptions in this section except **Remark 1** are from Chapter 1 of Alfonsi [71] with our notations. The CIR process $C = (C_t)_{t \geq 0}$ is the unique pathwise continuous solution to the Itô's SDE:

$$\mathrm{d}C_t = (a - rC_t)\mathrm{d}t + \sigma\sqrt{C_t}\mathrm{d}B_t, \quad t > 0 \tag{1}$$

subject to an initial condition $C_0 \geq 0$. Here, $a \geq 0$ is the source parameter, $r \geq 0$ is the reversion speed, and $\sigma \geq 0$ is the volatility. The SDE (1) has the affine drift term $a - rC_t$ and non-Lipschitz continuous diffusion term $\sigma\sqrt{C_t}$ (more rigorously, 1/2-Hölder continuous). This irregular setting of the diffusion term is key to modeling with the CIR process. In particular, the process $C$ is a.s. nonnegative at $t \geq 0$, and its positivity depends on the magnitude of volatility. If $\sigma^2 \leq 2a$ (low volatility), then the process is a.s. positive and hence does not reach the boundary value zero. By contrast, if $\sigma^2 > 2a$ (high volatility), then the process is not a.s. positive; it may reach the $t$ axis, hence becoming zero in a finite time. The volatility $\sigma$ size therefore determines the qualitative property of the CIR process.

Due to the particular form of the drift and diffusion terms in the SDE (1), the statistical moments of the CIR process can be obtained in a closed form. We explain that the CIR bridge—a time-dependent version of the process—has this tractability.

***Remark 1*** Sample paths of the high-volatility CIR process show on–off intermittency that is not observed for the low-volatility process (**Figure 1**). As discussed in **Section 4**, the CIR bridge identified from the data corresponds to the high-volatility case.

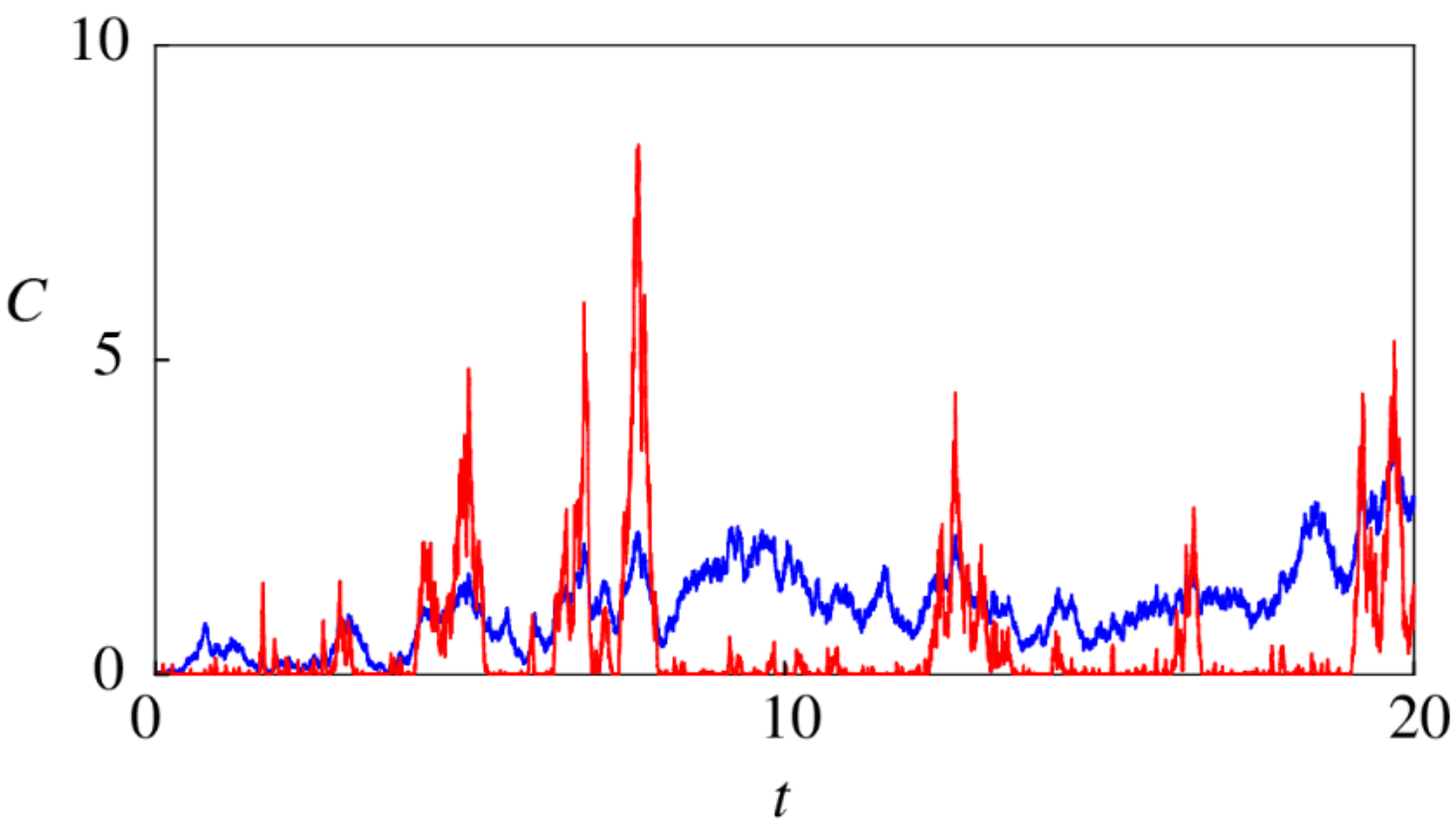

**Figure 1.** Comparisons of sample paths of CIR process starting from initial condition $C_0 = 0$ between low-volatility case (blue, $\sigma = 1$) and high-volatility case (red, $\sigma = 4$). We set $a = r = 1$. These sample paths are generated through the iVi scheme (**Section 4**).

### 2.2 CIR bridge

For each fixed terminal time $T > 0$, the proposed CIR bridge $X = (X)_{0 \le t \le T}$ is a pathwise solution to the following SDE:

$$\mathrm{d}X_t = \left(a(t) - h(t) X_t\right)\mathrm{d}t + \sigma\sqrt{h(t) X_t}\,\mathrm{d}B_t,\quad 0 < t < T \tag{2}$$

subject to the zero initial condition $X_0 = 0$. Here, $a : [0, T] \to [0, +\infty)$ is a bounded and continuous function, $h : [0, T) \to [0, +\infty)$ is a positive and continuous function that becomes unbounded at time $t = T$, and $(B_t)_{t \ge 0}$ is a 1-D standard Brownian motion. The SDE (2) is a time-inhomogeneous, localized version of (1). We computationally study a specific version of (2) in **Section 4**, whereas **Section 3** presents our theoretical investigation of the more generic form (2).

The design principle of the CIR bridge is to maintain the timescale consistency between the drift and diffusion terms, and it is analytically tractable as well as numerically computable. Specifically, our CIR bridge is a time-changed CIR process with a time-dependent source. Indeed, following the time change argument of Kazakevičius and Kononovicius [72], we can rewrite the SDE (2) as

$$\mathrm{d}X_t = \left(\frac{a(t)}{h(t)} - X_t\right)\mathrm{d}\left(H(t)\right) + \sigma\sqrt{X_t}\,\mathrm{d}B_{H(t)},\quad t > 0, \tag{3}$$

where $H$ is a primitive function of $h$.

**Table 1** compares the classical and present CIR bridges. The proposed bridge was motivated from the fact that the classical ones presented in the previous works [50, 63] are considered not appropriate for analysis of the fish migration that involves burst phenomenon, modelled by high-volatility cases, as reported in **Section 4.** Indeed, Example 2 in Chau et al. [63] imply that the classical bridge are valid for low-volatility cases (possibly "$q = \frac{2\kappa\mu}{\sigma^2} - 1$" needs to be nonnegative in the context of Chau et al. [63], where $\sigma$ is also the volatility in this reference). Theis finding motivated the modeling and analysis of a CIR bridge that applies to both high- and low-volatility cases. Moreover, the author did not find analytical formulae of the mean and variance of the classical CIR bridges. They are no longer affine processes, but highly nonlinear stochastic processes. By contrast, our CIR bridge is an affine process, with which the moments can be obtained analytically. In addition, how to apply the iVi scheme, which is a computational method tailored for affine processes, to the classical non-affine diffusion bridges are not trivial.

The increment in the right-hand side of (3) is related to the modified time $H(t)$ but not the time $t$, thus representing the time change with the common increment. The time change approach may be related to the scaling approach, where the drift and diffusion terms are scaled and divided by power functions of time [73]. The scaling approach is for designing diffusion processes with desired time-dependent probability densities, so it fundamentally differs from the time change approach. Nevertheless, both approaches have the same goal, which is to design a process by scaling variables.

From an ecological standpoint, the time change in (3) means that the migrating fishes have some biological clock so that they can sense environmental cues that govern the daytime duration, such as temperature and luminosity as some functions of the physical time. They are not directly considered in the proposed model but phenomenologically through the function $H$. Because an advantage of stochastic modeling is its conceptual nature, looking for a more detailed model for the fish migration would need a deeper reasoning that explicitly incorporates equations for environmental cues. This is beyond the scope of this paper but will be addressed in the future. Note that the proposed CIR bridge has not been studied so far, so it is an interesting mathematical object by itself.

***Remark 2*** In the context of martingale optimal transport and stochastic control, SDEs having state-dependent diffusion coefficients unbounded at the terminal time arise [74,75]. These models are different from ours, while both share the same principle that the unboundedness of diffusion can constrain the terminal condition in some sense.

***Remark 3*** Using a CIR process as a building block of the model for fish migration is motivated also from an ecological standpoint that the CIR and related processes arise as continuous-time limits of the population dynamics model based on the branching mechanism (e.g., Chapter 3 in Bansaye and Méléard [76]).

**Table 1.** Summary of comparison between the classical and present CIR bridges.

| | Classical [50][63] | Present |
|---|---|---|
| How was the model formulated? | By conditioning (h-transform) | By time change |
| Affine process? | No | Yes |
| Moments analytically available? | No | Yes (**Proposition 2**) |
| Considered for all volatility regimes? | No* | Yes |

* In [63].

## 3. Mathematical analysis

### 3.1 Well-posedness

We study the well-posedness of the SDE (2) and show that it admits a unique continuous pathwise nonnegative solution under certain conditions of coefficients $a$ and $h$ regardless of volatility $\sigma$. The following **Proposition 1** states that the SDE (2) is well-posed.

***Proposition 1***

*Assume that* $h:[0,T)\to+\infty$ *is continuously differentiable and satisfies*

$$\frac{\bar{h}}{T-t}\leq h(t)\leq\frac{\bar{h}}{T-t}+\omega,\quad 0<t<T \tag{4}$$

*with constants* $\bar{h}>0$ *and* $\omega\geq 0$. *Then, the SDE (2) admits a unique pathwise solution that is a.s. continuous in the time interval* $[0,T]$ *with the a.s. limit value* $\lim_{t\to T-0}X_t=0$.

In **Proposition 1**, the condition (4) is sufficient to achieve the terminal value. The following examples of $h$ satisfy it:

$$h(t)=\frac{c}{T-t}\ \text{ and }\ h(t)=\frac{1}{t+\varepsilon}+\frac{1}{T-t},\ \text{where}\ c,\varepsilon>0. \tag{5}$$

The primitives $H$ of the two cases are given by (up to a constant)

$$H(t)=c\ln\frac{1}{T-t}\ \text{ and }\ H(t)=\ln\frac{t+\varepsilon}{T-t}, \tag{6}$$

respectively. For the case of *P. altivelis* that we focus on, the blow up of the primitives $H$ near the terminal time $T$ means that the biological clock of migrating juveniles moves forward more rapidly near the sunset time so that the migration is terminated at or before that time.

A more complex form of $h$ may be theoretically possible as for bridges with additive diffusion terms [77]; nevertheless, our application (**Section 4**) shows that the models in (5) work reasonably. The form of $h$ assumed in (4) is satisfied by a limited class of functions. However, our results cover CIR bridges with unbounded drift and diffusion coefficients—a hitherto underexplored research area.

***Remark 4*** If $\varepsilon \to +0$ is assumed for the second equation in (5), then both the initial and terminal conditions are constrained. As for the terminal condition, this naturally leads to the initial value 0 of the bridge, which is assumed *a priori* in this paper. From a biological viewpoint, the condition $\varepsilon \to +0$ corresponds to the prohibition of migration near sunrise time, and this effect weakens as $\varepsilon$ increases.

***Remark 5*** From a purely theoretical standpoint, one may be interested in a CIR bridge with a positive terminal condition. For example, a naïve model of a CIR bridge ending at a prescribed value $X_T > 0$ is (with the initial condition $X_0 = 0$)

$$\begin{aligned} \mathrm{d}X_t &= \left(a(t) - h(t)(X_t - X_T)\right)\mathrm{d}t + \sigma\sqrt{h(t)X_t}\,\mathrm{d}B_t \\ &= \left\{\left(a(t) - h(t)X_t\right) + h(t)X_T\right\}\mathrm{d}t + \sigma\sqrt{h(t)X_t}\,\mathrm{d}B_t \end{aligned}, \quad 0 < t < T , \tag{7}$$

where the term $h(t)X_T \geq 0$ arises as an additional drift, and the model is an affine process that is well-posed at least in $[0,T)$. However, this model does not satisfy the continuity of $X$ at $T$ because the diffusion does not vanish at $T$ due to $X_T > 0$. In this view, our CIR bridge is tailored for the (biologically meaningful) case $X_T = 0$. This characteristic should be compared to bridges based on the h-transform that allow for $X_T > 0$ because of the nonlinear drift.

### 3.2 Moments and moment-generating function

As byproducts of **Proof of Proposition 1 in Appendix**, the closed-form formulae of the average and variance of the CIR bridge are obtained.

***Proposition 2***

*Under the assumption of* ***Proposition 1****, for any* $t \in [0,T]$*, the average* $\mathbb{E}[X_t]$ *and variance* $\mathbb{V}[X_t]$ *of the CIR bridge* $X$ *are given by*

$$\mathbb{E}[X_t] = \int_0^t a(s)\exp\left(-\int_s^t h(u)\,\mathrm{d}u\right)\mathrm{d}s \tag{8}$$

*and*

$$\mathbb{V}[X_t] = \sigma^2\int_0^t h(s)\mathbb{E}[X_s]\exp\left(-2\int_s^t h(u)\,\mathrm{d}u\right)\mathrm{d}s , \tag{9}$$

*respectively.*

We can also obtain the conditional moment-generating function in a closed form due to the forms of the drift and square-root diffusion.

***Proposition 3***

*Under the assumption of* ***Proposition 1****, for any* $\mu \in (0,T)$ *and* $\lambda < 2\sigma^{-2}$ *, the conditional moment-generating function of the CIR bridge*

$$\phi_\mu(t,x) = \mathbb{E}\left[\exp\left(\lambda X_{T-\mu}\right) \middle| X_t = x\right], \quad 0 \le t < T-\mu, \tag{10}$$

*is given by*

$$\phi(t,x) = \exp\left(\varphi_\mu(t)x + \psi_\mu(t)\right), \tag{11}$$

*where*

$$\psi_\mu(t) = \int_t^{T-\mu} a(u)\varphi_\mu(u)\mathrm{d}u \tag{12}$$

*and*

$$\varphi_\mu(t) = \frac{\lambda}{\left(1-\lambda\frac{\sigma^2}{2}\right)e^{\int_t^{T-\mu} h(u)\mathrm{d}u} + \lambda\frac{\sigma^2}{2}}. \tag{13}$$

*Moreover, it follows that*

$$\mathbb{E}\left[\exp\left(\lambda X_{T-\mu}\right)\right] = \exp\left(\psi_\mu(0)\right), \quad 0 \le t < T-\mu. \tag{14}$$

In **Proposition 3**, the assumption that $\lambda$ is small is necessary because the conditional moment-generating function will diverge as for the classical CIR case (e.g., Proposition 1.2.4 in Alfonsi [71]). If necessary, with **Proposition 3** in mind, one may be able to compute arbitrary moments and cumulants if necessary, by differentiating $\psi_\mu$ with respect to $\lambda$ required times and then taking the limit $\lambda \to +0$.

### 3.3 Numerical discretization

The CIR bridge is discretized using the iVi scheme developed by Abi Jaber [64] to deal with CIR processes with the full range of volatility values. This numerical scheme was designed by maintaining a consistency between the CIR process and its associated Riccati equation in a discrete sense. Th iVi scheme has been later extended by Abi Jaber and Attal [65] to compute a class of affine Volterra processes. Details of this scheme are available in Abi Jaber [64]. We present an algorithm that applies this scheme to the CIR bridge.

The time interval $[0,T]$ is discretized using the breaking points $t_i = i\Delta t$ ( $i = 0,1,2,\ldots,N$ ) with the time increment $\Delta t = \frac{T}{N}$ and some natural number $N \in \mathbb{N}$. The discretized $X$ at time $t_i$ is denoted as $X_{(i)}$. We compute each $X_{(i)}$ as follows, where $X_{(0)} = 0$; for $i = 0,1,2,\ldots,N-1$, we set

$$X_{(i+1)} = X_{(i)} + a_i\Delta t + b_i U_i + c_i Z_i, \tag{15}$$

where

$$a_i = a\left(t_i + \frac{1}{2}\Delta t\right), \quad b_i = -h\left(t_i + \frac{1}{2}\Delta t\right), \quad c_i = \sigma\sqrt{h\left(t_i + \frac{1}{2}\Delta t\right)}, \tag{16}$$

$$U_i \sim IG\left(\alpha_i, \left(\frac{\alpha_i}{\eta_i}\right)^2\right), \tag{17}$$

and

$$Z_i = \frac{1}{\eta_i}\left(U_i - \alpha_i\right). \tag{18}$$

Here, we set

$$\alpha_i = X_{(i)} \frac{e^{b_i \Delta t} - 1}{b_i} + \frac{a_i}{b_i}\left(\frac{e^{b_i \Delta t} - 1}{b_i} - \Delta t\right) \text{ and } \eta_i = c_i \frac{e^{b_i \Delta t} - 1}{b_i}, \tag{19}$$

and $IG(u,v)$ is the inverse Gaussian distribution with a probability density $\sqrt{\frac{v}{2\pi x^3}} e^{-\frac{v(x-u)^2}{2u^2 x}}$ ($x > 0$).

Discretization (15) is a time-inhomogeneous version of Algorithm 1 in Abi Jaber [64] and preserves the nonnegativity $X_{(i)} \geq 0$ through the direct application of Theorem 1.3 in Abi Jaber [64] regardless of the $\sigma$ value and the $h(t)$ size. This property is important for the current study because the net volatility $\sigma\sqrt{h(t)}$ becomes unbounded near the terminal time $T$ while the source $a(t)$ remains bounded in $[0,T]$. Consequently, the net volatility becomes arbitrarily larger than the source, rendering classical numerical methods, which assume the low-volatility condition, inapplicable. Another advantage of the iVi scheme is that the discretization method adopts fully explicit one-stage discretization per time step, as in the classical Euler–Maruyama scheme. Generating an inverse Gaussian random variable is not challenging; it can be achieved through acceptance rejection (Algorithm 2 [64]).

## 4. Application

### 4.1 Life history of *P. altivelis*

The diadromous fish *P. altivelis* typically has a one-year life history that includes migration between a river and a sea [32]. A generational change occurs in autumn (around October to November), when adults spawn in a downstream reach of a river and hatched larvae are transported to a sea (around November to February). The larvae grow up in the sea during winter, and juveniles migrate to a nearby river in late winter to early summer (around February to June). The juveniles grow in a midstream river reach from late spring to autumn (around May to October) and then migrate downstream for spawning. Among these life stages, this study focuses on the upstream migration of juveniles because its success (i.e., the number of migrants) directly influences population and growth dynamics, reproduction (summer to autumn), and spawning [78,79,80,81].

### 4.2 Study site and data

We study the upstream migration of *P. altivelis* in the Nagara River in Japan, a class A river in the Kiso River system that flows into the Atlantic Ocean. The *P. altivelis* in the Nagara River as a biological resource is a Globally Important Agricultural Heritage System due to its multifaced values [82].

Migrating *P. altivelis* has been counted for approximately 30 years at a fishway installed at the Nagara River Estuary Barrage (N 35.082, E 136.697; 5.4 km upstream from the Pacific Ocean. See **Figure 2**). Two data sets were obtained from the Japan Water Agency, Nagara River Estuary Barrage Management Office, Nagara River. The first data set contains *P. altivelis* juvenile count data every 10 min per day from February to June in 2023 and 2024. The second data set contains daily count data from February to June from 1995 to 2024. The fish counting by the Nagara River Estuary Barrage Management Office was conducted manually until 2019; it has been conducted through video recording enhanced by artificial intelligence since 2020 [83]. The daily count data are visualized on the homepage of the Nagara River Estuary Barrage Management Office without numerical data[2].

Fish counting during upstream migration has been conducted in many previous studies, including the Nagara River case, but most of them focus on daily data [84,85,86,87,88,89]; see also daily fish count data of five rivers summarized in Yoshioka [31]. Moreover, modeling using temporally fine time series data has not been addressed. Another motivation for studying the migration of *P. altivelis* in the Nagara River is the recent finding that 87.2% of the spawning fish in this river is natural [90] and hence many of the fish migrate between the river and the sea through the Nagara River Estuary Barrage.

**Figure 3** shows the daily fish count data from 1995 to 2024 (the total fish counts are presented in **Section A.3 in Appendix**). **Figure 3** suggests that the upstream migration at the study site starts in late February and finishes near the end of July, implying that the yearly migration duration is three to four months. **Figure 4** shows the 30-year average daily fish count, showing that both the average and standard deviation of the fish count peak in May. The coefficient of variation (CV) of each day varies between 2 and 5, showing that the fluctuation dominates the average in this data. Furthermore, the peak is related to the water temperature in the Nagara River as discussed below.

**Figure 5** shows the daily water temperature at the Oyabu Big Bridge Station (N 35.293, E 136.671. See **Figure 2**), 31 km upstream from the study site[3], which is the nearest observation point of water temperature of the Nagara River Estuary Barrage. Considering data availability, data from 2016 to 2024 are used in this study. According to **Figure 5**, significant fish migration (more than 1,000 fish per day) occurs at daily-average water temperatures of 11°C to 24°C, peaking at approximately 17°C to 18°C. Fish migration thus becomes significant at suitable water temperatures. The upstream and downstream migration of the thinlip gray mullet also shows this temperature dependence [91]. The study site and the Oyabu Big Bridge are expected to have different water temperatures because of their distance from each other. Nevertheless, the present analysis suggests that the water temperature at the Oyabu Big Bridge can be considered an environmental indicator for predicting the upstream migration of *P. altivelis*.

---

[2] https://www.water.go.jp/chubu/nagara/15_sojou/index.html (last accessed on May 14, 2025)

[3] http://www1.river.go.jp/cgi-bin/SiteInfo.exe?ID=405091285502120 (Last accessed on May 10, 2025)

Here, we present the 10 min fish count data, which are the main study target. **Figure 6** shows the daily 10 min fish count data in 2023 and 2024[4] on an ordinary logarithmic scale (see **Figure A1 in Appendix** for the plots in a nonlogarithmic scale), suggesting that fish migration occurs during daytime each day and becomes significant when the water temperature falls within the suitable range (**Figure 5**). Therefore, the time series of the 10 min fish count data per day can be considered a bridge, with the initial and terminal times being the sunrise and sunset times, respectively.

According to **Figure 6**, the fish count in every 10 min varies over four orders ranging from $O\left(10^0\right)$ to $O\left(10^4\right)$ with fluctuations. We fit the CIR bridge (with suitable normalization) to the sub-hourly fish migration data, as discussed in the next subsection. It turns out that the identified models are in the high-volatility conditions where the sample paths are intermittent. Moreover, the bursting nature of the 10 min data suggests that the manual counting of the fish migration, like that employed in the Nagara River before (e.g., counting fish manually for several minutes and several times), may give erroneous results because the empirical count would significantly depend on whether it is conducted during ON (bursting) or OFF (resting) states. Analyzing the intermittency of the CIR bridge is also important in this view.

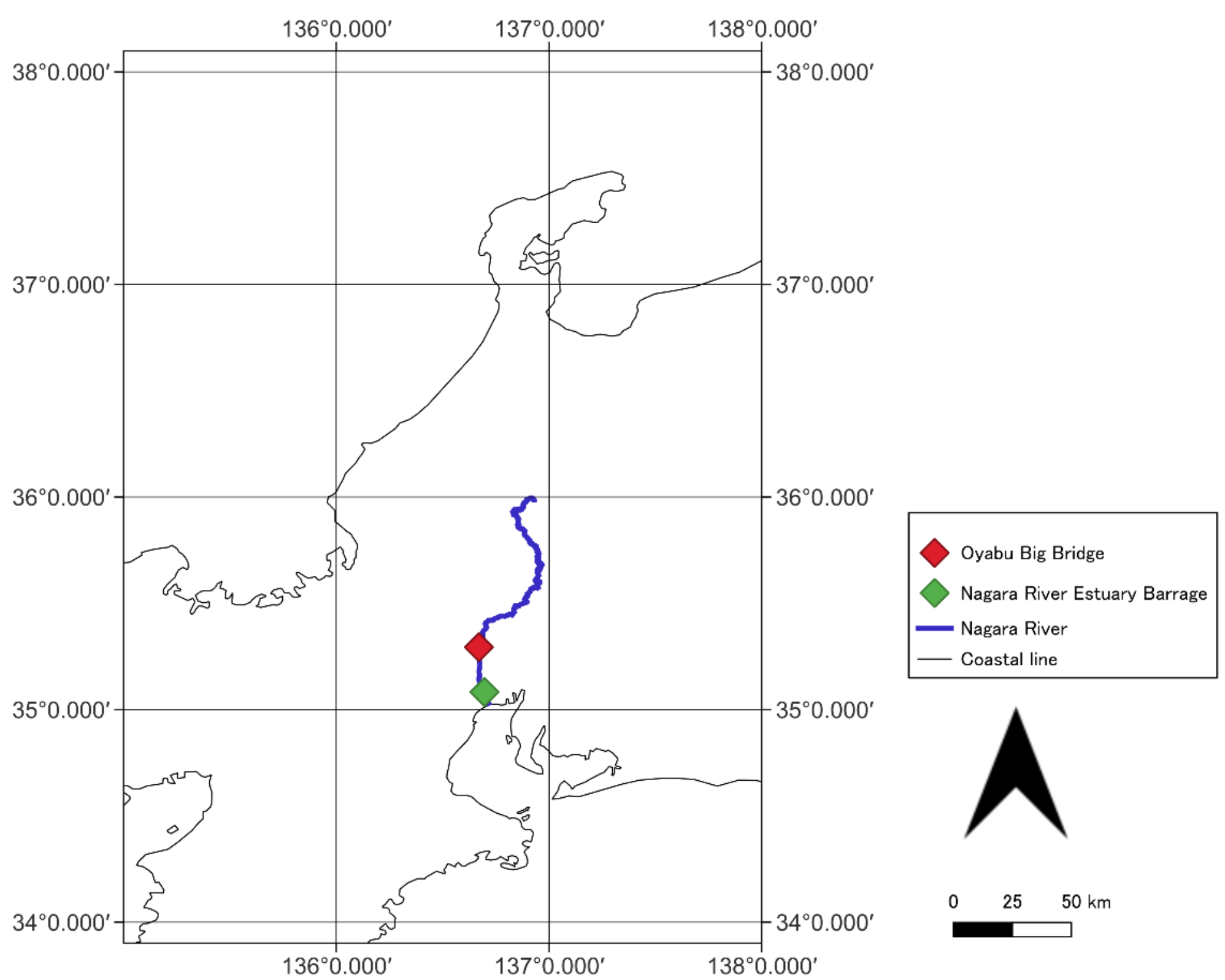

**Figure 2.** Map of the study area.

---

[4] The fish count on February 26 in 2023 was 0, but the observation window was not available in the data, so we left them blank in the panel (a) and excluded this day from the mathematical modeling in this paper.

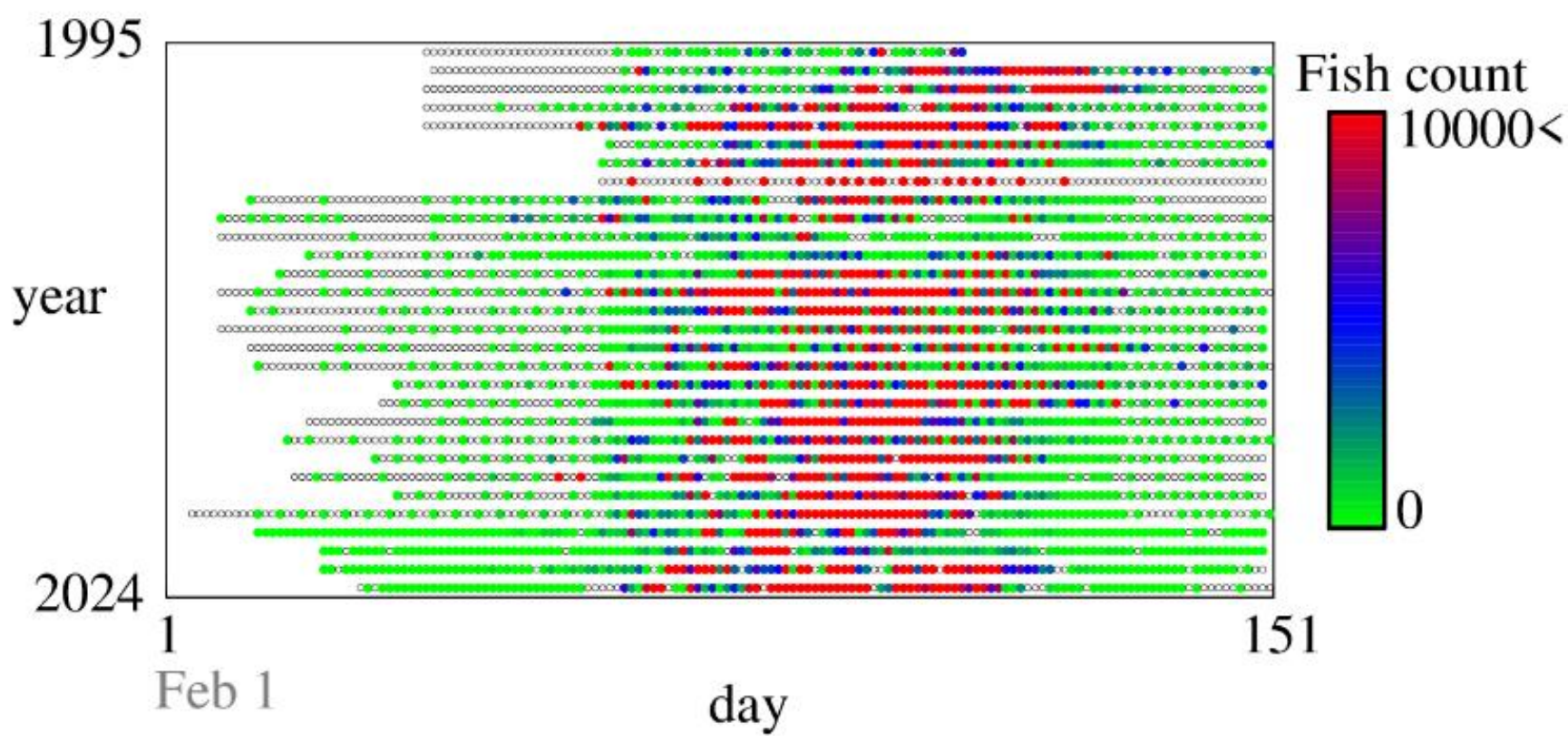

**Figure 3.** Daily fish count data from 1995 to 2024. Circles are filled if the data are available. Unfilled circles represent a fish count of 0.

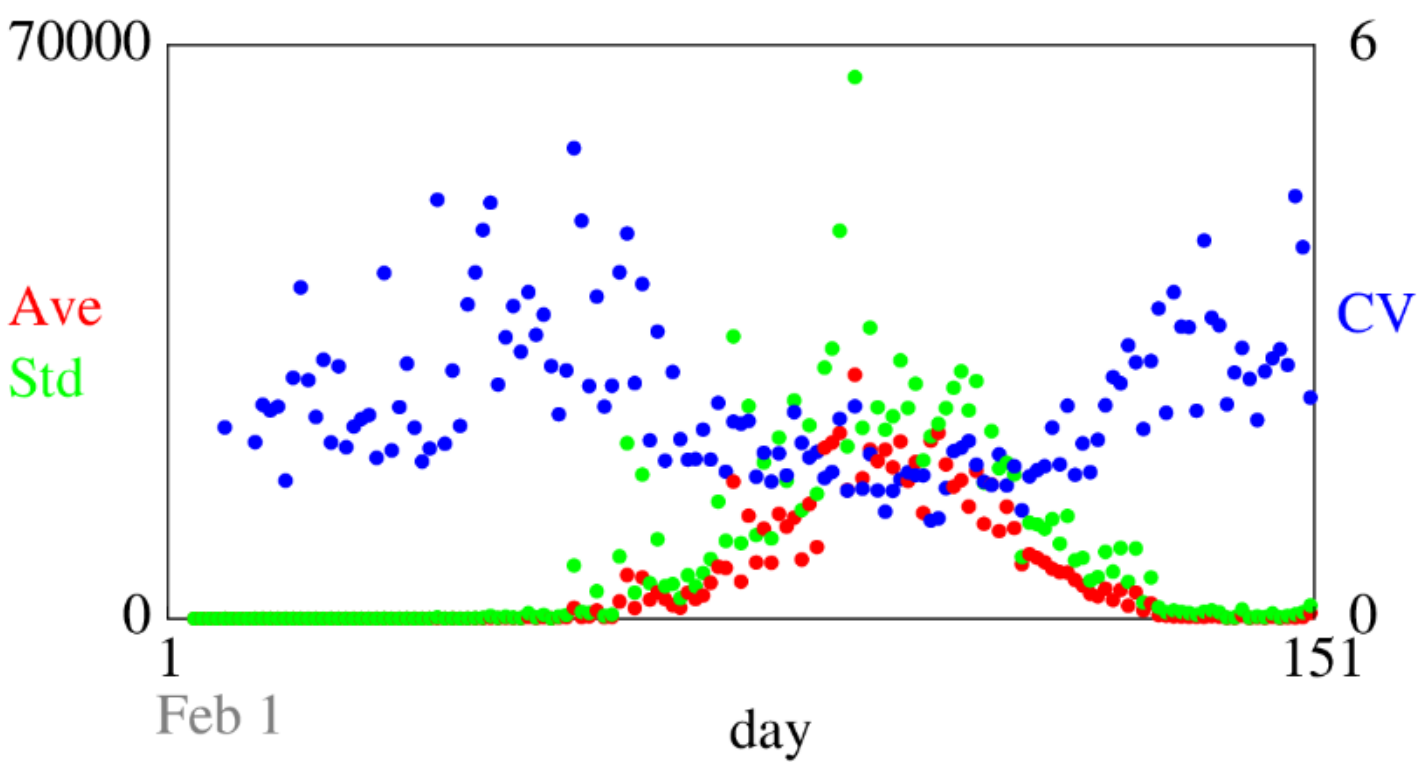

**Figure 4.** Average daily fish count data from 1995 to 2024. Ave: Average; Std: standard deviation; CV: coefficient of variation.

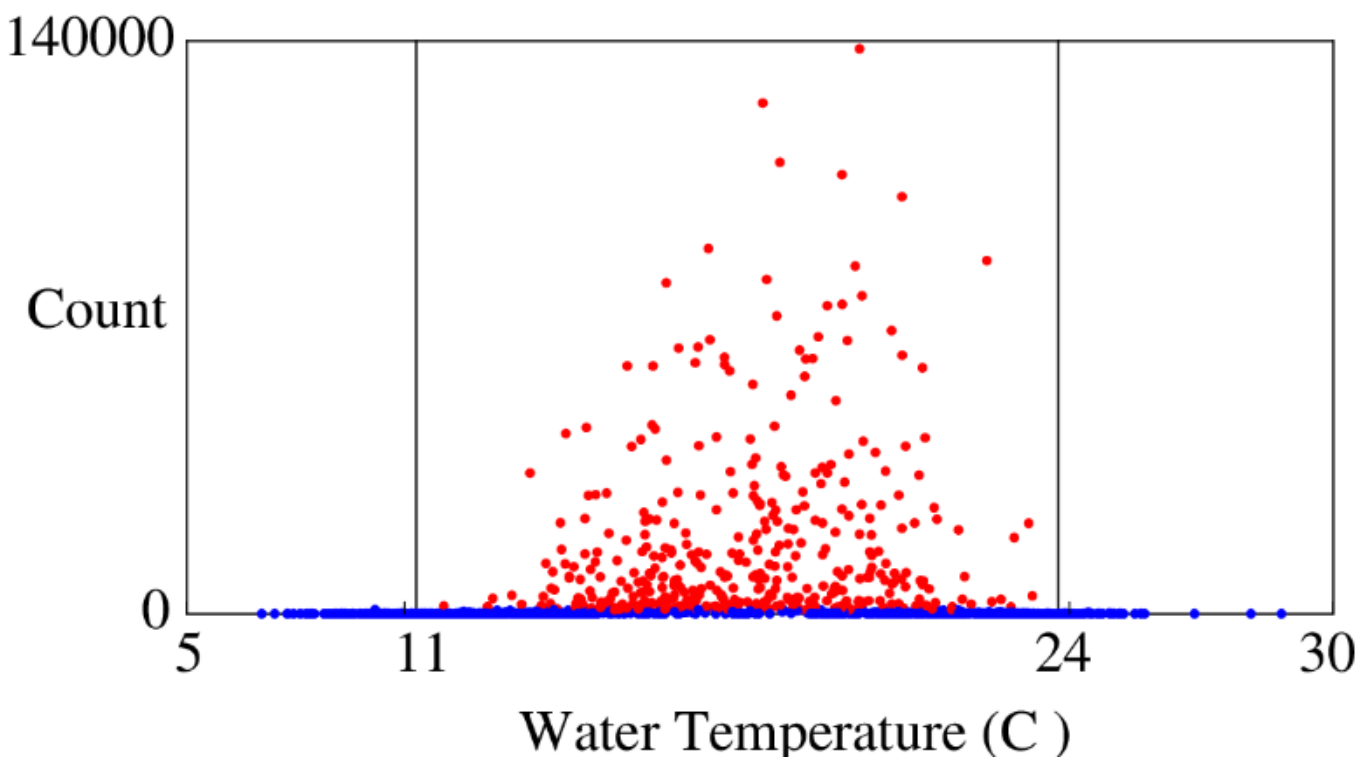

**Figure 5.** Relationship between daily-average water temperature and daily fish count. The red plots show counts exceeding 1,000 individuals.

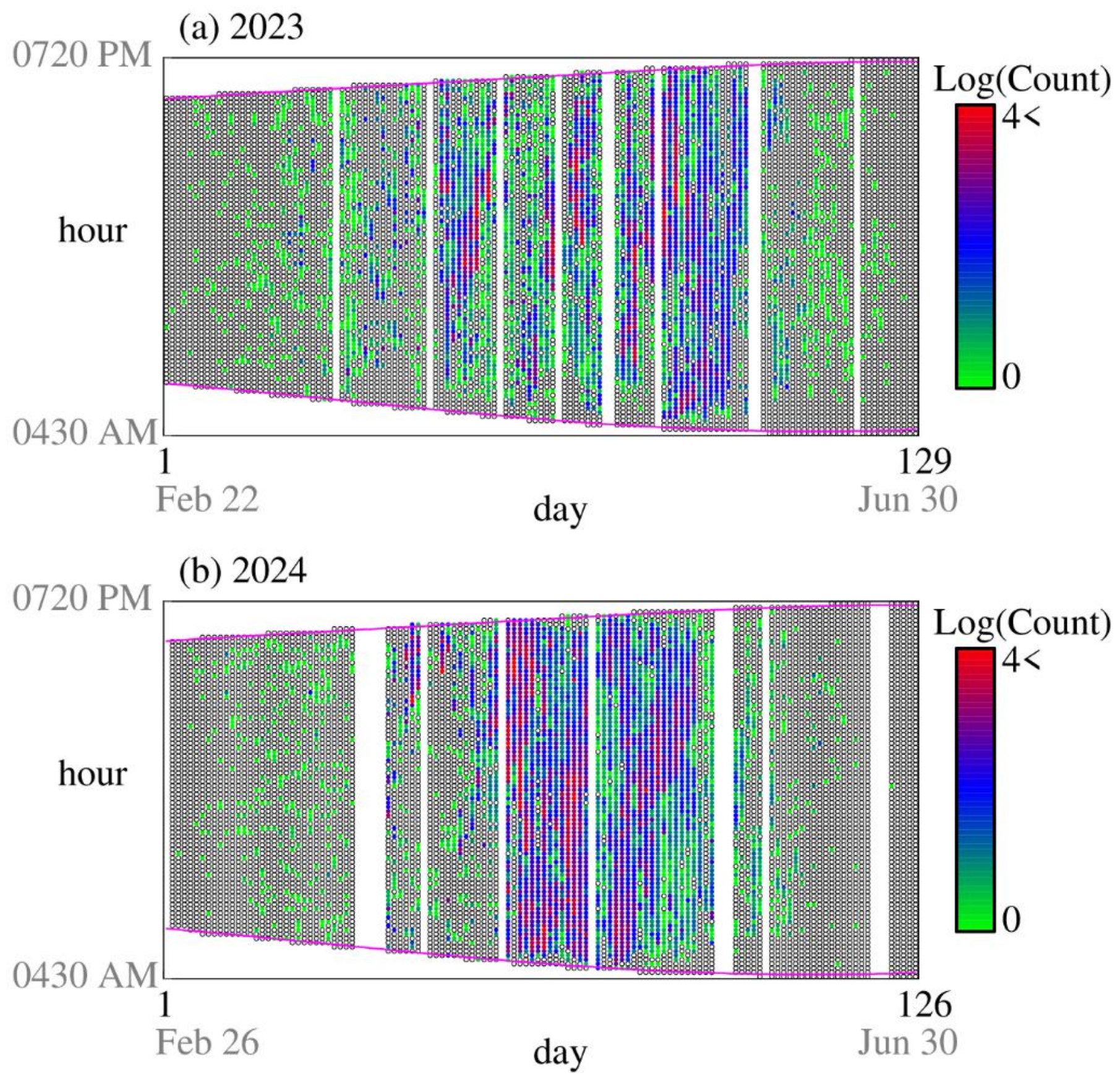

**Figure 6.** 10 min fish count data in (a) 2023 and (b) 2024 on an ordinary logarithmic scale. Circles are filled if the data are available. Unfilled circles represent a fish count of 0. The magenta curves represent sunrise and sunset times in Nagoya City, as obtained from the National Astronomical Observatory of Japan[5]. Plots on a nonlogarithmic scale are presented in **Section A.3**.

### 4.3 Parameter fitting

#### 4.3.1 CIR bridge formulation

The CIR bridge is estimated for each year (2023 or 2024). A year is fixed here. The sunrise and sunset times on the $k$ th day ( $k=1,2,3,...$ ) in the year are denoted as $T_{\mathrm{rise},k}$ and $T_{\mathrm{set},k}$ , respectively. We set $T_k = T_{\mathrm{set},k} - T_{\mathrm{rise},k} > 0$ . The unit-time fish count on day $k$ is assumed to follow the CIR bridge $X_{\cdot,k} = \left( X_{t,k} \right)_{T_{\mathrm{rise},k} \le t \le T_{\mathrm{set},k}}$ with the dimension of (Number/T):

$$\mathrm{d}X_{t,k} = \left( a_k - h_k\left(t\right) X_{t,k} \right) \mathrm{d}t + \sigma_k \sqrt{h_k\left(t\right) X_{t,k}}\, \mathrm{d}B_{t,k}\,,\quad T_{\mathrm{rise},k} < t < T_{\mathrm{set},k} \tag{20}$$

[5] https://eco.mtk.nao.ac.jp/koyomi/dni/ (Last accessed on May 15, 2025)

subject to the initial and terminal conditions $X_{T_{\mathrm{rise},k},k} = X_{T_{\mathrm{set},k},k} = 0$. Here, $\left(B_{t,k}\right)_{t\geq 0}$ ($\mathrm{T}^{1/2}$) is a 1-D standard Brownian motion that is mutually independent for different $k$, $a_k > 0$ (Number/$\mathrm{T}^2$) is a constant, a continuous function $h_k : \left[T_{\mathrm{rise},k}, T_{\mathrm{set},k}\right] \to [0,+\infty)$ (1/T), and $\sigma_k > 0$ ($\mathrm{Number}^{1/2}/\mathrm{T}^{1/2}$) is a volatility. The total daily fish count $S_k > 0$ (for now) [Number] on day $k$ is available in the data set. We set the nondimensional variable $Y_{s,k} = \dfrac{\Delta T}{S_k} X_{t,k}$, where $s = \dfrac{t - T_{\mathrm{rise},k}}{T_k} \in [0,1]$ and $\Delta T$ (T) is the time increment for normalization. By $X_{t,k} = \dfrac{S_k}{\Delta T} Y_{s,k}$, we transform (20) into the following SDE of $Y_{\cdot,k} = \left(Y_{s,k}\right)_{0\leq s\leq 1}$ for $0 < s < 1$:

$$\mathrm{d}\left(\frac{S_k}{\Delta T} Y_{s,k}\right) = \left(a_k - h_k(t)\frac{S_k}{\Delta T} Y_{s,k}\right)\mathrm{d}\left(T_k s\right) + \sigma_k \sqrt{h_k(t)\frac{S_k}{\Delta T} Y_{s,k}}\,\mathrm{d}B_{t,k} \tag{21}$$

or equivalently,

$$\mathrm{d}Y_{s,k} = \left(\frac{a_k T_k \Delta T}{S_k} - T_k h_k(t) Y_{s,k}\right)\mathrm{d}s + \sqrt{\frac{\Delta T}{S_k}}\sigma_k \sqrt{T_k h_k(t) Y_{s,k}}\,\mathrm{d}W_{s,k}, \tag{22}$$

subject to the initial and terminal conditions $Y_{0,k} = Y_{1,k} = 0$. Here, $\left(W_{t,k}\right)_{t\geq 0}$ is a 1-D standard Brownian motion that is mutually independent for different $k$. We use the following ansatz to obtain a common form of the SDE between different dates in the year; there exists nondimensional constants $a > 0$ and $\sigma > 0$ such that $a = \dfrac{a_k T_k \Delta T}{S_k}$ (-) and $\sigma = \sigma_k \sqrt{\dfrac{\Delta T}{S_k}}$ (-), and a nondimensional function $h : [0,1] \to [0,+\infty)$ such that $h(s) = T_k h_k(t)$ (-). Then, the SDE (22) reduces in the sense of law to a nondimensional form without $k$ dependence as follows:

$$\mathrm{d}Y_s = \left(a - h(s) Y_s\right)\mathrm{d}s + \sigma\sqrt{h(s) Y_s}\,\mathrm{d}W_s,\quad 0 < s < 1 \tag{23}$$

subject to the initial and terminal conditions $Y_0 = Y_1 = 0$. Here, $\left(W_t\right)_{t\geq 0}$ is a 1-D standard Brownian motion. The SDE (23) is (2) with $T = 1$. We can recover the parameters $a$ and $\sigma$ and the function $h$ for each day by calculating their $k$-dependent values because $T_{\mathrm{rise},k}$, $T_{\mathrm{set},k}$, and $S_k$ are available: $a_k = \dfrac{S_k}{T_k \Delta T} a$, $\sigma_k = \sqrt{\dfrac{S_k}{\Delta T}}\sigma$, and $h_k(t) = \dfrac{1}{T_k} h(s)$. We select $\Delta T = 10$ (min) because this study deals with 10 min fish count data. The choice of $\Delta T$ does not qualitatively affect the analysis results presented in the sequel.

The normalization employed here is controlled by the total daily fish count $S_k$ that varies several orders because $\Delta T$ is a constant and variations in $T_k$ are at most tens of percent. A larger $S_k$ implies larger $a_k$ and $\sigma_k$, and their behavior is qualitatively the same with those of $S_k$ and $S_k^{1/2}$, respectively, both being small outside moderate water temperatures as suggested from **Figure 5**.

#### 4.3.2 Fitted models

The full availability of average and variance as shown in **Proposition 2** motivates us to implement them into parameter fitting of the CIR bridge. We apply least-squares methods against the theoretical average and variance as explained below. An advantage of this fitting method is that parameter values can be identified without introducing any extra errors and approximations that are not avoidable in some other methods such as approximate maximum likelihood methods.

First, we identify the parameter $a$ and function $h$ using the average (8) through the least-squares method between the theoretical and empirical averages. This is because the average is independent of $\sigma$. Second, we identify the parameter $\sigma$ using the variance (9) through a least-squares method between the theoretical and empirical variances, where we assume the fitted results of $a$ and $h$. We assume the function $h$ of the forms given in (5) for $T=1$ because of the normalization. The CIR bridges corresponding to the left and right equations in (5) are called models 1 and 2, respectively. The $h$ functions for these models for $0<s<1$ are as follows:

$$h(t)=\frac{c}{1-s} \text{ (model 1) and } h(t)=\frac{1}{s+\varepsilon}+\frac{1}{1-s} \text{ (model 2), where } c,\varepsilon>0 . \tag{24}$$

The two models have the same total number of parameters and thus the same complexity. Moreover, these models admit closed-form expressions of the average and variance (**Section A.2**). The condition where diffusion dominates the source (high-volatility case) is expressed as follows through elementary calculations:

$$a<\frac{\sigma^2}{2}\min_{0\le s\le 1} h(s)=\frac{\sigma^2}{2}c \tag{25}$$

for model 1 and

$$a<\frac{\sigma^2}{2}\min_{0\le s\le 1} h(s)=\frac{4\sigma^2}{1+\varepsilon} \tag{26}$$

for model 2 if $\varepsilon\in(0,1)$ as in our case study discussed below.

### 4.4 Discussion

#### 4.4.1 Fitted results

**Tables 2 and 3** show the fitted parameter values of the CIR bridges. **Figures 7, 8, and 9** are comparisons of the theoretical and empirical averages, standard deviations, and CVs, respectively. **Table 4** gives the root mean squared errors (RMSEs) between the empirical and theoretical results calculated by the models, demonstrating that model 2 outperforms model 1 except in the fitting of the average in 2024. Thus, in the rest of **Section 4**, "the model" or "CIR bridge" means "model 2."

**Figures 7 and 8** show that both the average and variance are close to 0 near the initial and terminal times, suggesting that fish migration can be modeled using a bridge. Furthermore, the fitted results in the figures are reasonably spread around the empirical curves. The fluctuations of the empirical data are considered due to the smallness of the number of sample paths (120 to 130 data sample paths each year). Empirical CVs have an average of approximately 3 and vary between 2 and 6 (**Figure 9**), exhibiting qualitatively the same tendency as the daily count data. This tendency is captured by the CIR bridge with the fitted parameter values in both years. According to the fitted model, the average and variance peaks occur at the normalized time $s > 1/2$, which corresponds to the afternoon each day. These peaks are closer to the terminal time 1, corresponding to sunset time, in 2024 than in 2023.

The fitted values of $a$ and $\sigma$ do not considerably differ between 2023 and 2024, with a relative difference of approximately 10%, whereas those of $\varepsilon$ in 2023 are over four times that in 2024. Thus, the mean reversion in 2024 is stronger than that in 2023, particularly near the initial time 0 (corresponding to sunrise time). This difference appears as the smaller CV in 2023 (**Figure 9**). According to (25) and (26), all fitted models correspond to the high-volatility case, where diffusion (always) dominates the source, leading to intermittent sample paths (See **Figure 10** in the next subsection). The sub-hourly fish migration at the study site is therefore an intermittent phenomenon corresponding to the high-volatility case. This is a mathematical characterization of the fish migration phenomenon found by using the CIR bridge.

We also discuss environmental fluctuations that potentially affect the upstream migration of *P. altivelis*. Spagnolo et al. [69] and Valenti et al. [92] suggested that environmental noises may critically affect qualitative behavior such as stability of ecosystems both in space and time. This paper studied only the temporal dynamics of fish migration at a fixed observation point, and hence no spatial information is contained in the collected data. At the present time, it is technically difficult to observe spatio-temporal fish migration dynamics because it needs installation of many video cameras along a river with a fine spatial interval. Nevertheless, our data analysis could give a hint for better understanding the fish migration in future. As shown in Spagnolo et al. [69], time series at a fixed point of a spatio-temporal ecosystem model qualitatively depends on the intensity and spatial correlation of background environmental noise, which in our context would be water temperature, water quality, and discharge. We infer that, if the water temperature is dominant in the background noise and it drives the spatio-temporal migration dynamics of *P. altivelis*, then measuring its spatial correlation would be a way toward building an advanced mathematical model of fish migration. Indeed, our analysis suggested that water temperature is a key environmental factor for the migration of *P. altivelis*. We consider that surveying water temperature along the river will be technically more realistic than directly investigating spatio-temporal fish migration dynamics.

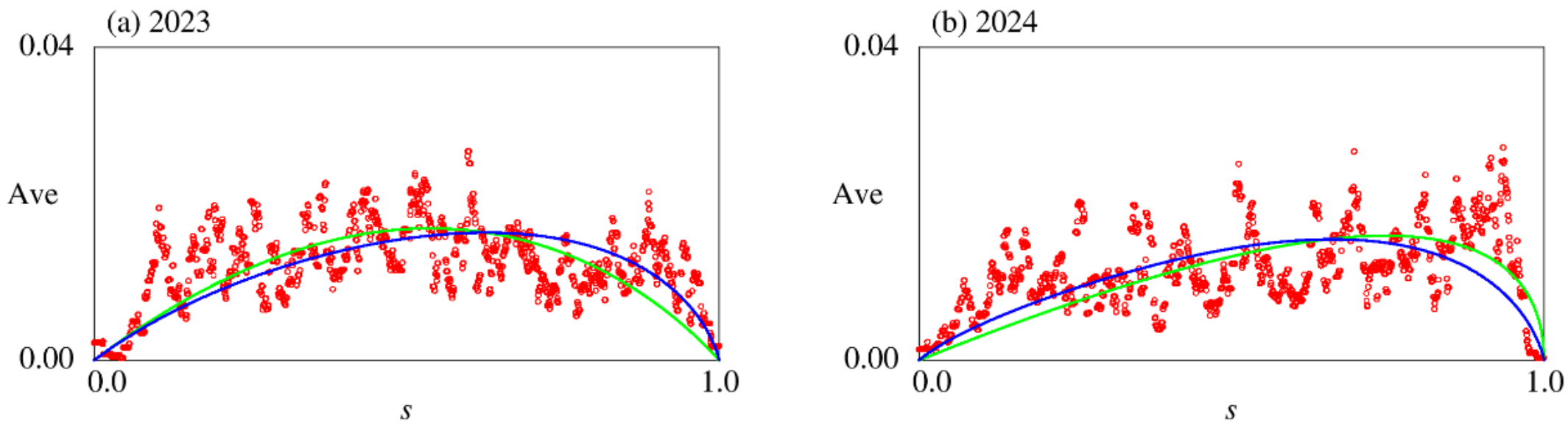

**Figure 7.** Comparison of theoretical (curves: green for 2023, blue for 2024) and empirical (circles) averages (Ave) of $Y$ in (a) 2023 and (b) 2024.

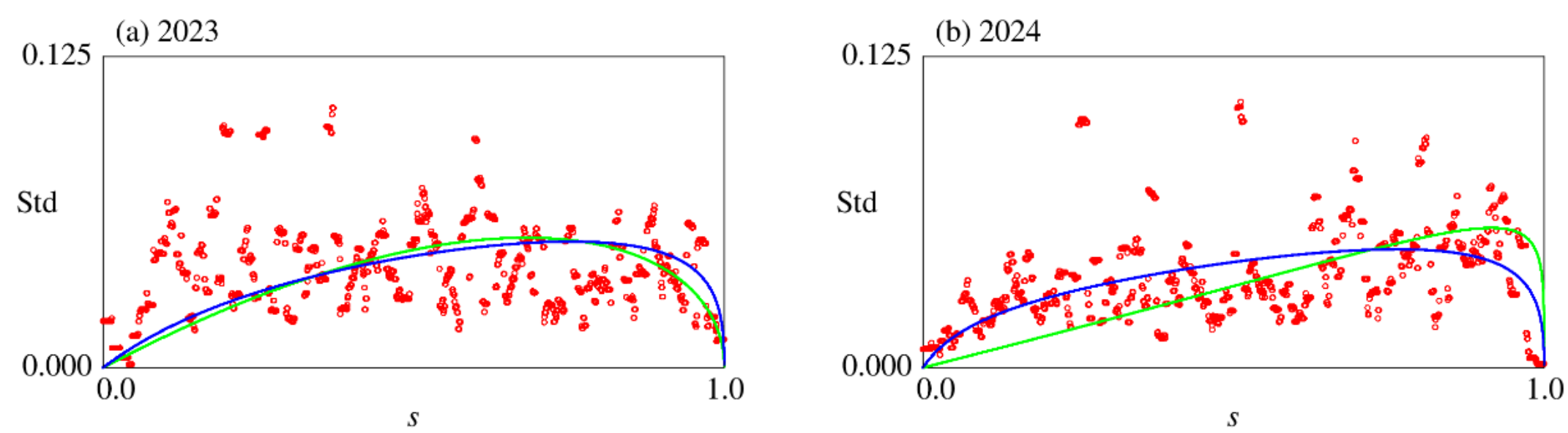

**Figure 8.** Comparison of theoretical (curves: green for model 1, blue for model 2) and empirical (circles) standard deviations (Std) of $Y$ in (a) 2023 and (b) 2024.

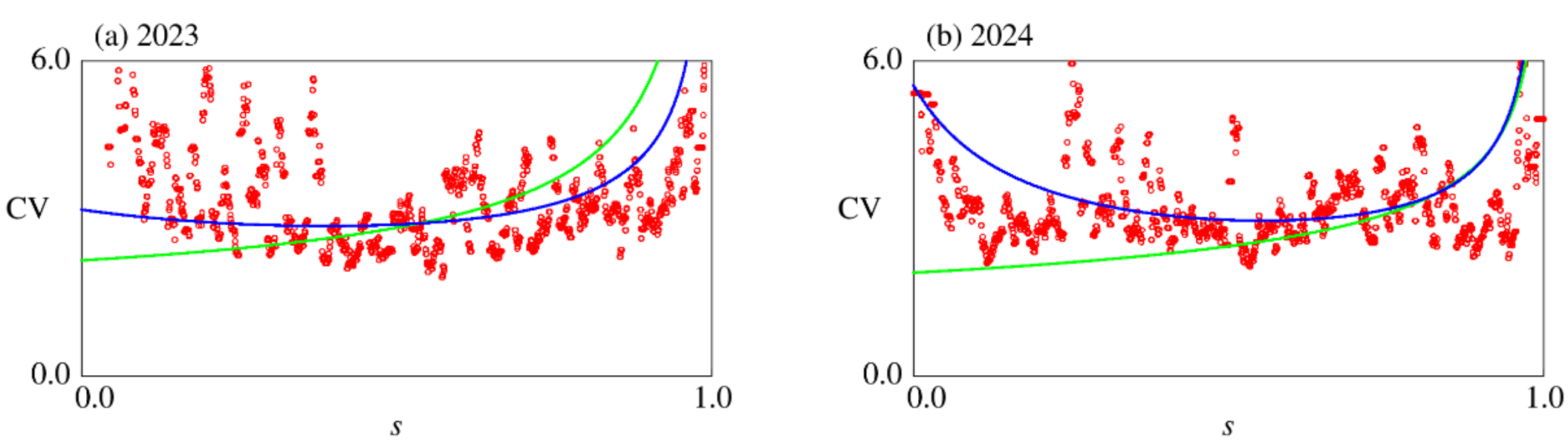

**Figure 9.** Comparison of theoretical (curves: green for 2023, blue for 2024) and empirical (circles) coefficient of variations (CV) of $Y$ in (a) 2023 and (b) 2024.

**Table 2.** Fitted parameter values for model 1.

| | 2023 | 2024 |
|---|---|---|
| $a$ | 0.0600 | 0.0322 |
| $c$ | 1.6473 | 0.5137 |
| $\sigma$ | 0.5942 | 0.6962 |

**Table 3.** Fitted parameter values for model 2.

| | 2023 | 2024 |
|---|---|---|
| $a$ | 0.0599 | 0.0666 |
| $\varepsilon$ | 0.3837 | 0.0808 |
| $\sigma$ | 0.5775 | 0.5523 |

**Table 4.** RMSEs of models 1 and 2.

| | Model 1 (2023) | Model 1 (2024) | Model 2 (2023) | Model 2 (2024) |
|---|---|---|---|---|
| Average | 4.329.E-03 | 4.807.E-03 | 4.283.E-03 | 4.870.E-03 |
| Standard deviation | 2.168.E-02 | 2.132.E-02 | 2.089.E-02 | 1.763.E-02 |

### 4.4.2 Convergence analysis

We study convergence of the numerical method against the CIR bridge in 2023. **Figure 10** shows sample paths of the CIR bridge computed via the iVi scheme with $\Delta t = 0.001$. They show some intermittency in the sample paths, with several bursts observed during each path. In this paper, pseudo-random numbers are generated by the Mersenne twister [93].

The time-dependent average and standard deviation of the model are correctly captured by the iVi scheme with 1,000,000 sample paths (**Figure 11**). For an analysis of the convergence of the iVi scheme against different computational conditions (sample size and time step), **Tables 5 and 6** summarize the maximum errors between the theoretical and computed averages and standard deviations. The iVi scheme reproduces the theoretical average and standard deviation if the sample size is sufficiently large. An excessively small time step increases the error, as in the results for the 10,000 samples. However, an overly large sample size also increases the error; they therefore must be balanced. The errors in the average are one order of magnitude smaller than those in the standard deviation. The convergence speed of the iVi scheme is difficult to quantify, but the diagonal elements in **Tables 5 and 6** suggest that increasing both the sample size and time step by at least ten times halves the error. **Tables A2 and A3 in Section A.3 of Appendix** show the time-average error for the cases corresponding to **Tables 5 and 6**, respectively, showing that the time-average error is several times to one order smaller than the maximum one.

We also explore how the terminal condition 0 is approximated in the numerical solutions. **Table 7** shows the error of the terminal values measured using the theoretical average value 0 and the computational one at the terminal time 1. The error decreases as the resolution in time increases and again indicates the need to balance the sample size.

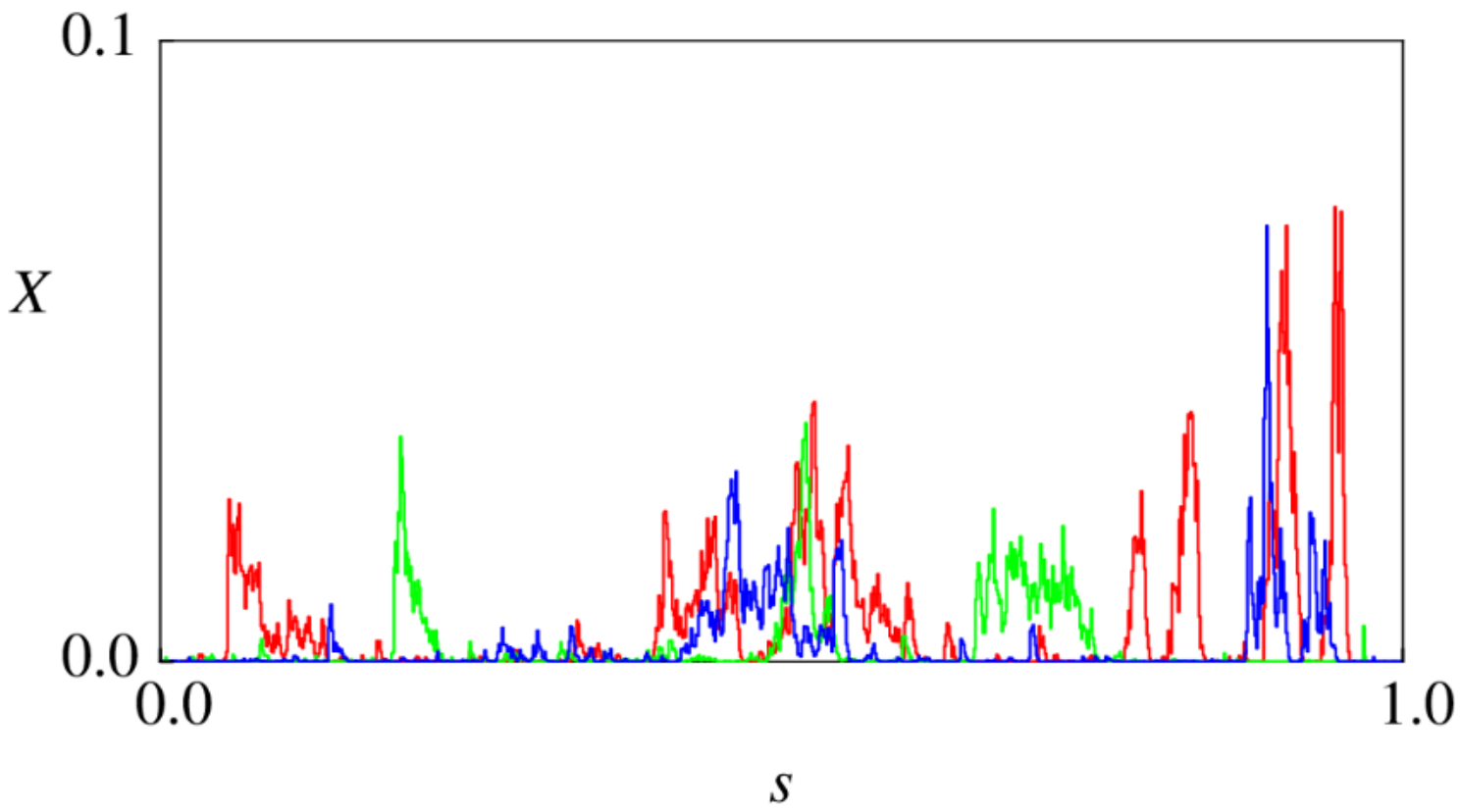

**Figure 10.** Sample paths of CIR bridge with parameter values in 2023. Distinct colors represent different sample paths.

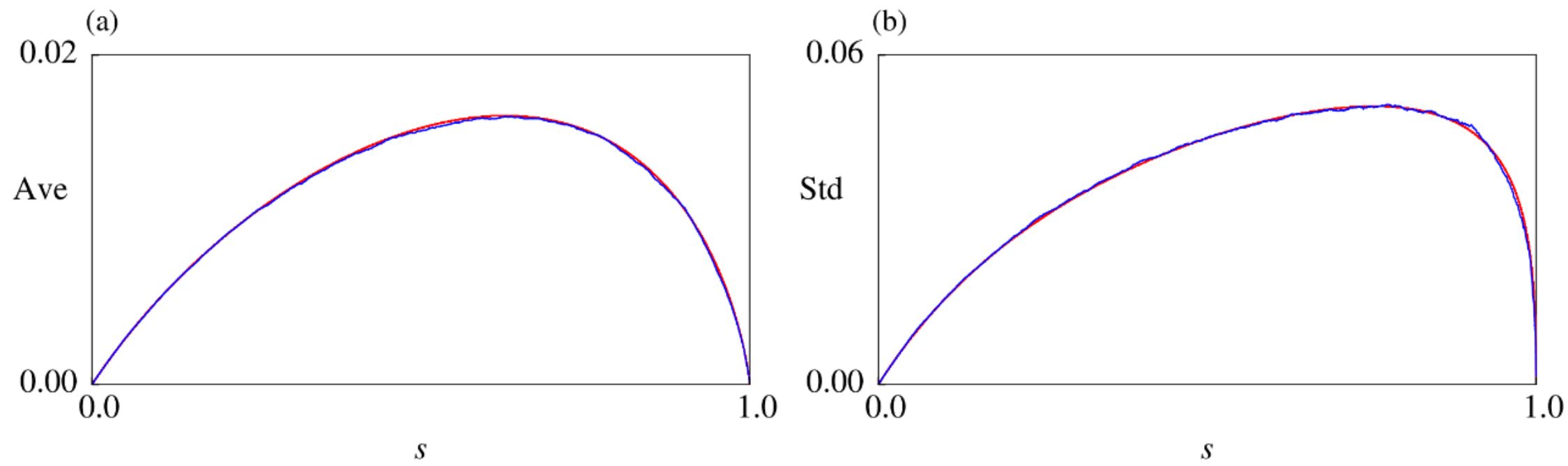

**Figure 11.** Comparison of theoretical (red) and computed (blue) results for (a) average and (b) standard deviation.

**Table 5.** Maximum error between computed and theoretical averages across all time steps.

| | | Sample size | | |
|---|---|---|---|---|
| | | 10,000 | 100,000 | 1,000,000 |
| Time step | 100 | 9.347.E-04 | 5.816.E-04 | 5.894.E-04 |
| | 1000 | 5.997.E-04 | 2.776.E-04 | 8.099.E-05 |
| | 10000 | 9.043.E-04 | 1.880.E-04 | 6.617.E-05 |

**Table 6.** Maximum error between computed and theoretical standard deviations across all time steps.

| | | Sample size | | |
|---|---|---|---|---|
| | | 10,000 | 100,000 | 1,000,000 |
| Time step | 100 | 7.206.E-03 | 4.587.E-03 | 4.811.E-03 |
| | 1000 | 2.731.E-03 | 2.808.E-03 | 1.935.E-03 |
| | 10,000 | 4.824.E-03 | 2.268.E-03 | 1.155.E-03 |

**Table 7.** Error of terminal values measured using theoretical terminal value 0 and computational average at terminal time 1.

| | | Sample size | | |
|---|---|---|---|---|
| | | 10,000 | 100,000 | 1,000,000 |
| Time step | 100 | 5.698.E-04 | 5.816.E-04 | 5.894.E-04 |
| | 1,000 | 4.119.E-05 | 7.157.E-05 | 8.099.E-05 |
| | 10,000 | 7.976.E-06 | 8.376.E-06 | 9.268.E-06 |

#### 4.4.3 Burst analysis

We also investigate the pathwise properties of the CIR bridge. Given the intermittency of its computed sample paths (**Figure 10**), we evaluate the number of bursts per sample path using a previously reported method [94,95]; see related models in Phillips et al. [96]. In this method, a burst is an event where the value of the state variable exceeds some prescribed threshold. This threshold height of burst events is denoted as $\underline{X} > 0$. We adopt the threshold duration $\underline{T} > 0$, which is sufficiently larger than $\Delta t$ but sufficiently smaller than 1, meaning a bust event is additionally required to continue not shorter than $\underline{T}$. We introduce $\underline{T}$ because a preliminary computation implies that small $X$ fluctuations near the threshold $\underline{X}$ spuriously increase the count of burst events. $\Delta t = 0.0001$, $\underline{T} = 0.02$, and $\underline{X} = 0.01$ unless otherwise specified. This setting targets burst events that persist longer than 10–20 min with heights at least comparable to the average of $X$. We use 100000 independent paths for sampling. The selected time increment $\Delta t$ and the number of sample paths are chosen based on the results presented in **Tables 5 through 7** where the terminal condition is satisfied with relative maximum errors at the order of $O\left(10^{-6}\right)$ and average and standard deviations are several %. The time-average errors are much smaller than these results as shown in **Tables A2 and A3 in Section A.3 of Appendix**.

**Table 8** shows the cumulants, CVs, and the maximum and minimum total numbers of burst events per day for 2023 and 2024. **Figure 12** shows the computed probability distributions (PDs) of the total number of burst events. We also examine cases where the threshold duration $\underline{T}$ for detecting burst events is doubled to $\underline{T} = 0.04$ and the threshold height $\underline{X} = 0.02$. In all four cases, the computed probability density functions (PDFs) of the burst events exhibit unimodality and positive skewness, and the CV is approximately 1. **Figure 12** indicates that 0, 1, or 2 bursts are likely to occur each day. In an extreme case, the computational results show that 6 to 9 burst events per day are possible, consistent with **Figure 6**, which depicts some sparseness in fish counts. Thus, the proposed CIR bridge is considered be able to handle both sparse- and dense-burst cases.

**Tables 9** shows the computed cumulants, CVs, and the maximum and minimum durations of burst events. **Figure 13** shows the computed PDFs of burst event durations. In all cases, the average burst duration is approximately 0.1, which corresponds to about one hour; the CV is again around 1. The CIR bridge for 2023 predicts a larger number and longer durations of burst events on average than the 2024 bridge (**Tables 8 and 9**). A similar tendency applies to the variance. According to the proposed mathematical model, the fish migration in 2023 is a bit more frequent and erratic than that in 2024 in view of burst events, but their difference is at most about 10% in terms of the statistics. Fish migration in these years is therefore not considered to be critically different in view of the proposed mathematical model. **Figure 13** shows that the modulation of the threshold values $\underline{T}$ and $\underline{X}$ examined here does not qualitatively affect the shape of the computed PDFs of the burst duration; particularly, the influences of $\underline{X}$ is at most few percents in view of the statistics presented in **Table 9**. This implies certain robustness of the method in the earlier studies [94,95].

Finally, the analysis of burst events in the CIR bridge suggested that this model, and possibly the fish migration of *P. altivelis* at least at the study site as well, has intermittency and presents burst events several times in one day. Manually counting the migrating fish would be better designed if this intermittency is accounted for. Of course, such an approach should depend on monetary and labor costs. Optimizing these costs and benefits of accurately counting the fish will be an important research topic in future.

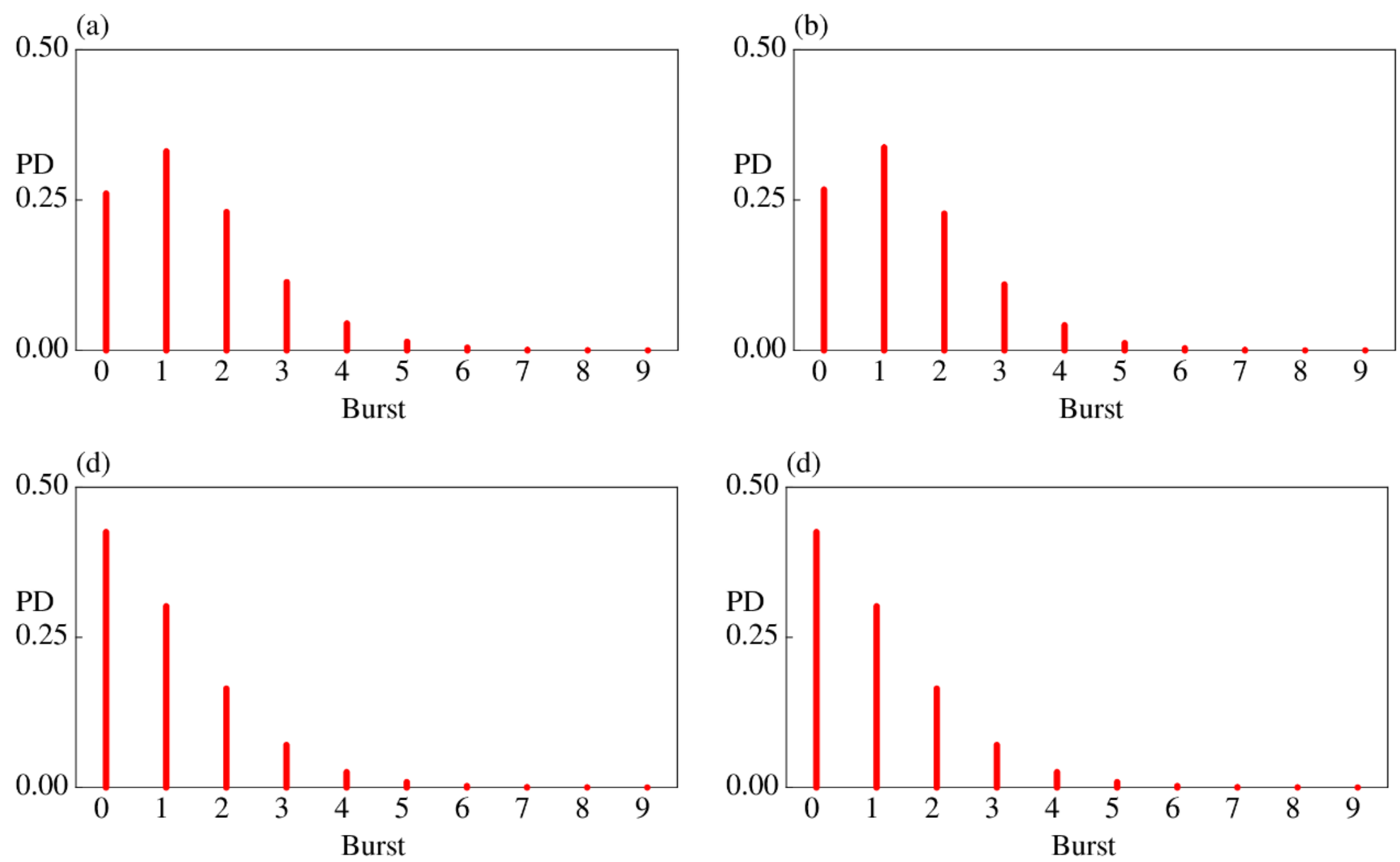

**Figure 12.** PDs of total number of bursts of CIR bridge for (a) 2023, (b) 2024, (c) 2024 ( $\underline{T}$ doubled), and (d) 2024 ( $\underline{X}$ doubled).

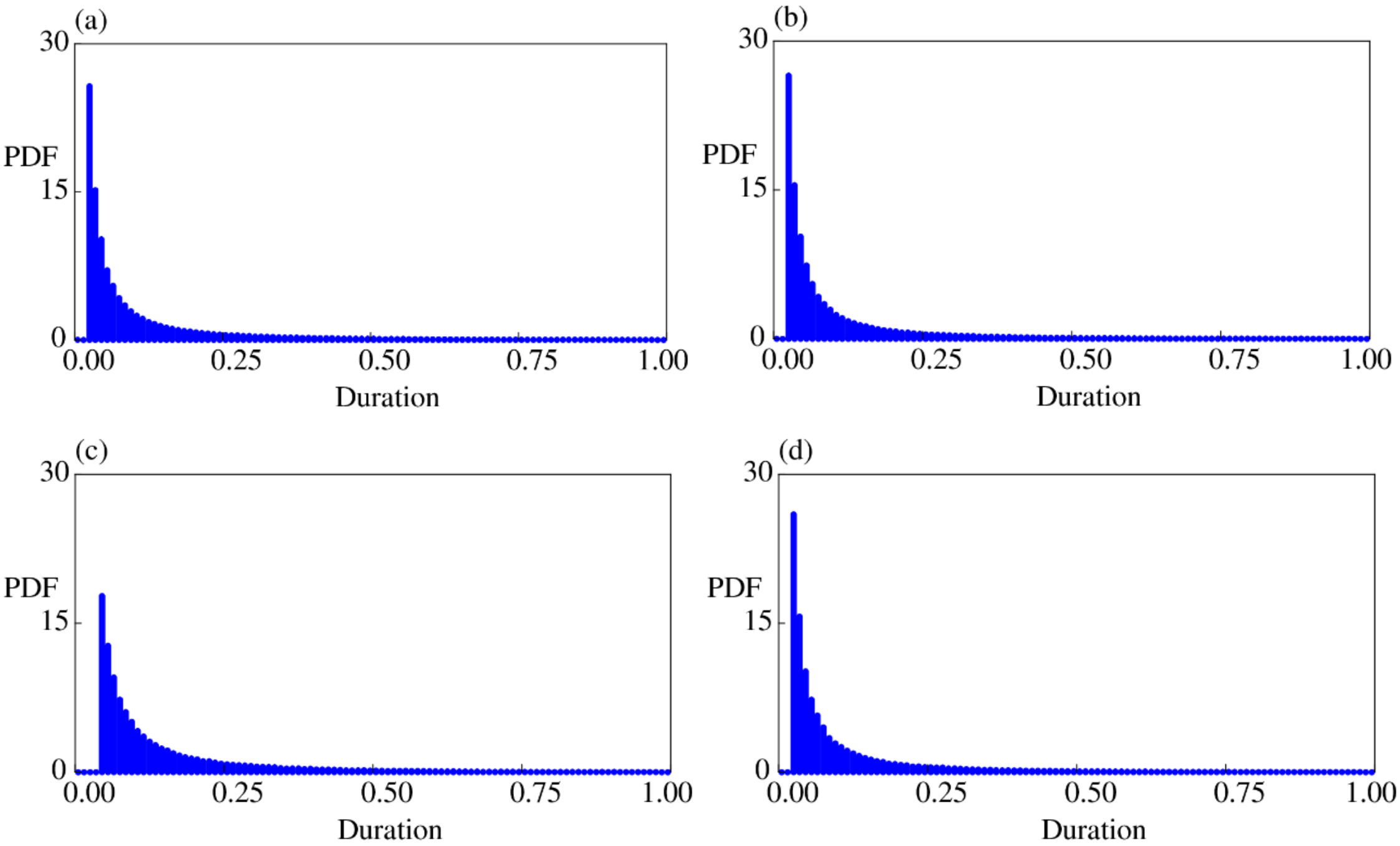

**Figure 13.** PDFs of burst durations of CIR bridge for (a) 2023, (b) 2024, (c) 2024 ( $\underline{T}$ doubled), and (d) 2024 ( $\underline{X}$ doubled).

**Table 8.** Computed cumulants, CVs, and maximum and minimum total numbers of burst events per day.

| Statistics | 2023 | 2024 | 2024 ($\underline{T}$ doubled) | 2024 ($\underline{X}$ doubled) |
|---|---|---|---|---|
| Average | 1.419.E+00 | 1.377.E+00 | 7.993.E-01 | 1.011.E+00 |
| Variance | 1.564.E+00 | 1.487.E+00 | 7.702.E-01 | 1.341.E+00 |
| Skewness | 9.369.E-01 | 9.256.E-01 | 1.061.E+00 | 1.320.E+00 |
| Kurtosis | 9.108.E-01 | 8.673.E-01 | 1.037.E+00 | 1.975.E+00 |
| CV | 8.815.E-01 | 8.856.E-01 | 1.098.E+00 | 1.146.E+00 |
| Maximum | 9 | 9 | 6 | 9 |
| Minimum | 0 | 0 | 0 | 0 |

**Table 9.** Computed cumulants, CVs, and maximum and minimum durations of burst events.

| Statistics | 2023 | 2024 | 2024 ($\underline{T}$ doubled) | 2024 ($\underline{X}$ doubled) |
|---|---|---|---|---|
| Average | 9.230.E-02 | 8.671.E-02 | 1.290.E-01 | 8.623.E-02 |
| Variance | 1.114.E-01 | 1.021.E-01 | 1.169.E-01 | 9.965.E-02 |
| Skewness | 2.974.E+00 | 3.056.E+00 | 2.466.E+00 | 3.026.E+00 |
| Kurtosis | 1.081.E+01 | 1.180.E+01 | 7.403.E+00 | 1.167.E+01 |
| CV | 1.207.E+00 | 1.177.E+00 | 9.060.E-01 | 1.156.E+00 |
| Maximum | 9.715.E-01 | 9.771.E-01 | 9.771.E-01 | 9.756.E-01 |
| Minimum | 2.010.E-02 | 2.010.E-02 | 4.000.E-02 | 2.010.E-02 |

## 5. Conclusion

We proposed and investigated the CIR bridge as the solution to an SDE through its application to fish migration. The CIR bridge is well-posed and has a closed-form average and variance. Its parameter values were fully determined in the application study, which focused on sub-hourly *P. altivelis* count data in 2023 and 2024 where the initial and terminal times of the SDE were determined to be the sunrise and sunset times, respectively. The daily count data demonstrated the temperature dependence of the daily fish count, giving additional information on the fish migration phenomenon. We found that the fish migration on the sub-hourly scale in the study site could be characterized as an intermittent phenomenon corresponding to a high-volatility case. We computationally validated the iVi scheme and used it to quantify such intermittency, focusing on burst events in sample paths. In summary, in this paper we found that the CIR bridge as a well-posed mathematical model can be applied to the upstream migration of *P. altivelis* as a highly stochastic and intermittent biological phenomenon.

The proposed mathematical approach would apply to other diadromous fishes with available fine fish-count time-series data. The affine property of the CIR bridge makes it ideal for planning problems with fishery resources, such as fish harvesting problems at a fixed point in a river. A limitation of the proposed bridge is that its terminal time is fixed. This assumption is suited to our setting, but it is not in the bridge with a randomized terminal time [97]. The applicability of the CIR bridge is not limited to fish migration; it can be applied to many other real-world phenomena where the target variable is nonnegative and positive in only a part of the whole domain. For example, the bridge can be used to model discharge in an intermittent river environment [98]. Applications of the CIR bridge and its related models to ecological and environmental phenomena are currently being undertaken by the author. CIR-type processes have been found to be able to describe the noise-induced stability in stochastic systems [99]. In future, the proposed CIR bridge can be incorporated into stability analysis of some predator-prey model where the predator is *P. altivelis* and prey is its feed resources, such as benthic algae and aquatic insects. Further collecting the sub-hourly fish migration data at the study site is also an important task, with which we will be able to deeper understand this interesting biological phenomenon.

## Appendix

### A.1 Proofs

***Proof of Proposition 1***

The coefficient $a$ is a bounded and continuous function in $[0,T]$. Therefore, through a regularization argument for the square-root diffusion term analogous to that in Proposition 1.2.1 of Alfonsi [71], the SDE (2) admits a unique continuous pathwise solution in any time interval $[0,T-\delta]$ with any $\delta \in (0,T)$, during which $h$ is bounded and smooth. Therefore, the problem is to study the solution at the terminal time $t = T$. Throughout this proof, $C > 0$ is a sufficiently large constant that depends only on $a$, $\bar{h}$, $\omega$, $T$, and its value will be updated line by line when necessary.

In the time interval $[0,T-\delta]$ with any $\delta \in (0,T)$, we obtain the average $\mathbb{E}[X_t]$ using the ordinary differential equation obtained by directly taking the expectation with $\mathbb{E}[X_0] = 0$:

$$\frac{\mathrm{d}}{\mathrm{d}t}\mathbb{E}[X_t] = a(t) - h(t)\mathbb{E}[X_t], \tag{27}$$

from which we obtain

$$\mathbb{E}[X_t] = \int_0^t a(s)\exp\left(-\int_s^t h(u)\mathrm{d}u\right)\mathrm{d}s. \tag{28}$$

Similarly, the variance $\mathbb{V}[X_t] = \mathbb{E}\left[\left(X_t - \mathbb{E}[X_t]\right)^2\right]$ is determined using another ordinary differential equation (analogous to Eq. (14) in Mostafa and Allen [100]) with $\mathbb{V}[X_0] = 0$:

$$\frac{\mathrm{d}}{\mathrm{d}t}\mathbb{V}[X_t] = \sigma^2\mathbb{E}[X_t] - 2h(t)\mathbb{V}[X_t], \tag{29}$$

from which we obtain

$$\mathbb{V}[X_t] = \sigma^2\int_0^t h(s)\mathbb{E}[X_s]\exp\left(-2\int_s^t h(u)\mathrm{d}u\right)\mathrm{d}s. \tag{30}$$

We investigate formulae (28) and (30) near the terminal time $T$ and set $\bar{a} = \max_{0 \le s \le T} a(s)$. For the average (28), for $t \in (0,T)$, we have

$$
\begin{aligned}
&0 \le \mathbb{E}\left[X_t\right] \\
&\le \bar{a}\int_0^t \exp\left(-\int_s^t h(u)\,\mathrm{d}u\right)\mathrm{d}s \\
&= \bar{a}\exp\left(-\int_0^t h(u)\,\mathrm{d}u\right)\int_0^t \exp\left(\int_0^s h(u)\,\mathrm{d}u\right)\mathrm{d}s \\
&\le \bar{a}\exp\left(-\int_0^t \frac{\bar{h}}{T-u}\,\mathrm{d}u\right)\int_0^t \exp\left(\int_0^s\left(\omega+\frac{\bar{h}}{T-u}\right)\mathrm{d}u\right)\mathrm{d}s\,. \\
&\le C\exp\left(-\int_0^t \frac{\bar{h}}{T-u}\,\mathrm{d}u\right)\int_0^t \exp\left(\int_0^s \frac{\bar{h}}{T-u}\,\mathrm{d}u\right)\mathrm{d}s \\
&= C\exp\left(-\bar{h}\ln\left(\frac{T}{T-t}\right)\right)\int_0^t \exp\left(\bar{h}\ln\left(\frac{T}{T-s}\right)\right)\mathrm{d}s \\
&= C\left(\frac{T-t}{T}\right)^{\bar{h}}\int_0^t\left(\frac{T}{T-s}\right)^{\bar{h}}\mathrm{d}s
\end{aligned} \tag{31}
$$

We have

$$
\begin{aligned}
\int_0^t\left(\frac{T}{T-s}\right)^{\bar{h}}\mathrm{d}s &= T^{\bar{h}}\int_0^t (T-s)^{-\bar{h}}\,\mathrm{d}s \\
&= \frac{T^{\bar{h}}}{1-\bar{h}}\left(T^{1-\bar{h}}-(T-t)^{1-\bar{h}}\right)\ \left(\text{if } \bar{h}\ne 1\right). \\
&= T\ln\left(\frac{T}{T-t}\right)\ \left(\text{if } \bar{h}=1\right)
\end{aligned} \tag{32}
$$

By (31) and (32), we obtain

$$
0 \le \lim_{t\to T-0}\mathbb{E}\left[X_t\right] \le \lim_{t\to T-0} C\left(\frac{T-t}{T}\right)^{\bar{h}}\int_0^t\left(\frac{T}{T-s}\right)^{\bar{h}}\mathrm{d}s = 0\,, \tag{33}
$$

and hence $\lim_{t\to T-0}\mathbb{E}\left[X_t\right]=0$. For the variance (30), for $t\in(0,T)$, we have

$$
\begin{aligned}
&0 \le \mathbb{V}\left[X_t\right] \\
&= \sigma^2\exp\left(-2\int_0^t h(u)\,\mathrm{d}u\right)\int_0^t h(s)\,\mathbb{E}\left[X_s\right]\exp\left(2\int_0^s h(u)\,\mathrm{d}u\right)\mathrm{d}s \\
&\le \sigma^2\exp\left(-2\int_0^t \frac{\bar{h}}{T-u}\,\mathrm{d}u\right)\int_0^t\left(\omega+\frac{\bar{h}}{T-s}\right)\mathbb{E}\left[X_s\right]\exp\left(\int_0^s\left(2\omega+\frac{2\bar{h}}{T-u}\right)\mathrm{d}u\right)\mathrm{d}s \\
&\le \sigma^2\exp\left(-2\int_0^t \frac{\bar{h}}{T-u}\,\mathrm{d}u\right)\int_0^t\left(\frac{T}{T-s}\omega+\frac{\bar{h}}{T-s}\right)\mathbb{E}\left[X_s\right]\exp\left(2\omega T+2\int_0^s \frac{\bar{h}}{T-u}\,\mathrm{d}u\right)\mathrm{d}s\,. \\
&\le C\exp\left(-2\int_0^t \frac{\bar{h}}{T-u}\,\mathrm{d}u\right)\int_0^t \frac{1}{T-s}\mathbb{E}\left[X_s\right]\exp\left(2\int_0^s \frac{\bar{h}}{T-u}\,\mathrm{d}u\right)\mathrm{d}s \\
&\le C\left(\frac{T-t}{T}\right)^{2\bar{h}}\int_0^t\left(\frac{T}{T-s}\right)^{2\bar{h}+1}\mathbb{E}\left[X_s\right]\mathrm{d}s
\end{aligned} \tag{34}
$$

If $\bar{h}=1$, then by (33), we obtain

$$
\begin{aligned}
\int_0^t \left(\frac{T}{T-s}\right)^{2\bar{h}+1} \mathbb{E}[X_s]\mathrm{d}s &\le \int_0^t \left(\frac{T}{T-s}\right)^{3} C\left(\frac{T-s}{T}\right) T \ln\left(\frac{T}{T-s}\right)\mathrm{d}s \\
&\le C\int_0^t \left(\frac{1}{T-s}\right)^{2} \ln\left(\frac{1}{T-s}\right)\mathrm{d}s \\
&= O\left(\frac{1}{T-t}\ln\left(\frac{1}{T-t}\right)\right)
\end{aligned} , \tag{35}
$$

where the last line holds true if $t$ is close to $T$. By (34) and (35), we obtain $\lim_{t\to T-0}\mathbb{V}[X_t]=0$ if $\bar{h}=1$.

If $\bar{h}\neq 1$, then by (32) we have

$$
\begin{aligned}
\int_0^t \left(\frac{T}{T-s}\right)^{2\bar{h}+1} \mathbb{E}[X_s]\mathrm{d}s &\le \int_0^t \left(\frac{T}{T-s}\right)^{2\bar{h}+1} C\left(\frac{T-s}{T}\right)^{\bar{h}} \left(T^{1-\bar{h}}-(T-s)^{1-\bar{h}}\right)\mathrm{d}s \\
&\le C\int_0^t (T-s)^{-\bar{h}-1}\frac{T^{1-\bar{h}}-(T-s)^{1-\bar{h}}}{1-\bar{h}}\mathrm{d}s \\
&= \frac{C}{1-\bar{h}}\left(T^{1-\bar{h}}\frac{(T-s)^{-\bar{h}}-T^{-\bar{h}}}{\bar{h}}-\int_0^t (T-s)^{-2\bar{h}}\mathrm{d}s\right)
\end{aligned} . \tag{36}
$$

The last integral in (36) is evaluated as

$$
\int_0^t (T-s)^{-2\bar{h}}\mathrm{d}s = \begin{cases} O\left(\ln\left(\dfrac{1}{T-t}\right)\right) & (2\bar{h}=1) \\ O\left((T-t)^{1-2\bar{h}}\right) & (2\bar{h}\neq 1) \end{cases} . \tag{37}
$$

Therefore, by (34) and (36)-(37), we obtain $\lim_{t\to T-0}\mathbb{V}[X_t]=0$ if $\bar{h}\neq 1$.

The variation of the constant formula shows that the unique pathwise solution of SDE (2) satisfies

$$
X_t = \int_0^t a(s)\exp\left(-\int_s^t h(u)\mathrm{d}u\right)\mathrm{d}s + \int_0^t \exp\left(-\int_s^t h(u)\mathrm{d}u\right)\sigma\sqrt{h(s)X_s}\,\mathrm{d}B_s \tag{38}
$$

in the time interval $[0,T-\delta]$ with any $\delta\in(0,T)$, which is continuous in this interval as well. By (28), we can rewrite (38) as

$$
X_t - \mathbb{E}[X_t] = \int_0^t \exp\left(-\int_s^t h(u)\mathrm{d}u\right)\sigma\sqrt{h(s)X_s}\,\mathrm{d}B_s \, (=M_t), \tag{39}
$$

where $M=(M_t)_{0\le t\le T-\delta}$ is a martingale. According to (30), $M$ is a uniformly square-integrable martingale in the time interval $[0,T-\delta]$ with any $\delta\in(0,T)$ because for any $t\in[0,T-\delta]$,

$$
\mathbb{E}[M_t]=0 \text{ and } \mathbb{E}\left[(M_t)^2\right]=\mathbb{E}\left[(X_t-\mathbb{E}[X_t])^2\right]=\mathbb{V}[X_t]<\bar{V}<+\infty \tag{40}
$$

with a constant $\bar{V}>0$ independent of $\delta$. Moreover, we have

$$
\lim_{t\to T-0}\mathbb{E}[M_t]=0 \text{ and } \lim_{t\to T-0}\mathbb{E}\left[(M_t)^2\right]=0 . \tag{41}
$$

This, along with the martingale convergence theorem (e.g., Proof of Proposition 2 in Yoshioka and Yamazaki [54]), there a.s. exists a random variable $M_T = \lim_{t\to T-0} M_t$, which is actually 0 due to (41). With

(38) in mind, this observation proves that the SDE (2) admits a unique pathwise nonnegative solution that is a.s. continuous on the time interval $[0,T]$ with the a.s. limit $\lim_{t\to T-0} X_t = 0$.

□

***Proof of Proposition 3***

The CIR bridge is a time-inhomogeneous affine process (e.g., Filipović [101]). Therefore, for any sufficiently smooth function $\phi:[0,T-\mu]\times[0,+\infty)\to\mathbb{R}$ with $\mu\in(0,T)$, the infinitesimal generator $\mathbb{L}\phi$ is given by (Theorem 2.13 in Filipović [101])

$$\mathbb{L}\phi(t,x)=\frac{\partial\phi(t,x)}{\partial t}+\left(a(t)-h(t)x\right)\frac{\partial\phi(t,x)}{\partial x}+\frac{\sigma^2}{2}h(t)x\frac{\partial^2\phi(t,x)}{\partial x^2}. \tag{42}$$

The conditional characteristic function then satisfies

$$\mathbb{L}\phi_\mu(t,x)=0 \quad \text{for} \quad 0\le t<T-\mu \quad \text{and} \quad x>0 \tag{43}$$

subject to the following terminal condition with a small $\lambda$ whose range is determined later in the proof:

$$\phi_\mu(T-\mu,x)=\exp(\lambda x) \quad \text{for} \quad x>0. \tag{44}$$

We guess the solution of the form $\phi_\mu(t,x)=\exp\left(\varphi_\mu(t)x+\psi_\mu(t)\right)$ with $\varphi_\mu,\psi_\mu:[0,T-\mu]\to\mathbb{R}$ such that $\varphi_\mu(T-\mu)=\lambda$ and $\psi_\mu(T-\mu)=0$. Substituting this $\phi_\mu$ into (43) yields

$$\mathbb{L}\phi_\mu(t,x)=\phi_\mu(t,x)\left\{\left(\frac{\mathrm{d}\varphi_\mu(t)}{\mathrm{d}t}-h(t)\varphi_\mu(t)+\frac{\sigma^2}{2}h(t)\left\{\varphi_\mu(t)\right\}^2\right)x+\frac{\mathrm{d}\psi_\mu(t)}{\mathrm{d}t}+a(t)\varphi_\mu(t)\right\}, \tag{45}$$

and hence

$$\frac{\mathrm{d}\varphi_\mu(t)}{\mathrm{d}t}-h(t)\varphi_\mu(t)+\frac{\sigma^2}{2}h(t)\left\{\varphi_\mu(t)\right\}^2=0 \quad \text{for} \quad 0\le t<T-\mu \tag{46}$$

and

$$\frac{\mathrm{d}\psi_\mu(t)}{\mathrm{d}t}+a(t)\varphi_\mu(t)=0 \quad \text{for} \quad 0\le t<T-\mu. \tag{47}$$

These time-backward equations can be written in a time-forward form by introducing $t=T-s$, $\bar{\varphi}_\mu(s)=\varphi_\mu(T-s)$, and $\bar{\psi}_\mu(s)=\psi_\mu(T-s)$:

$$\frac{\mathrm{d}\bar{\varphi}_\mu(s)}{\mathrm{d}s}=h(T-s)\left\{-\bar{\varphi}_\mu(s)+\frac{\sigma^2}{2}\left(\bar{\varphi}_\mu(s)\right)^2\right\} \quad \text{for} \quad \mu<s\le T \tag{48}$$

and

$$\frac{\mathrm{d}\bar{\psi}_\mu(s)}{\mathrm{d}s}=a(T-s)\bar{\varphi}_\mu(s) \quad \text{for} \quad \mu<s\le T \tag{49}$$

with $\bar{\varphi}_\mu(\mu)=\lambda$ and $\bar{\psi}_\mu(\mu)=0$.

Applying the transformation of variables $\bar{\varphi}_\mu(s) = \dfrac{1}{\gamma(s)}$ to (48) yields the linear ordinary differential equation

$$\frac{\mathrm{d}\gamma(s)}{\mathrm{d}s} = h(T-s)\left\{\gamma(s) - \frac{\sigma^2}{2}\right\} \quad \text{for} \quad \mu < s \le T \tag{50}$$

with $\gamma(\mu) = \lambda^{-1}$. Assuming $\lambda^{-1} > \dfrac{\sigma^2}{2}$ yields

$$\begin{aligned}
\gamma(s) &= \lambda^{-1} e^{\int_\mu^s h(T-u)\mathrm{d}u} - \frac{\sigma^2}{2}\int_\mu^s e^{\int_v^s h(T-u)\mathrm{d}u} h(T-v)\,\mathrm{d}v \\
&= e^{\int_\mu^s h(T-u)\mathrm{d}u}\left\{\lambda^{-1} + \frac{\sigma^2}{2}\left[e^{-\int_\mu^v h(T-u)\mathrm{d}u}\right]_\mu^s\right\} \\
&= e^{\int_\mu^s h(T-u)\mathrm{d}u}\left\{\lambda^{-1} + \frac{\sigma^2}{2}\left(e^{-\int_\mu^s h(T-u)\mathrm{d}u} - 1\right)\right\} \qquad \text{for} \quad \mu < s \le T, \\
&= \left(\lambda^{-1} - \frac{\sigma^2}{2}\right) e^{\int_\mu^s h(T-u)\mathrm{d}u} + \frac{\sigma^2}{2} \\
&> 0
\end{aligned} \tag{51}$$

and hence

$$\bar{\varphi}_\mu(s) = \varphi_\mu(T-s) = \frac{\lambda}{\left(1 - \lambda\dfrac{\sigma^2}{2}\right) e^{\int_\mu^s h(T-u)\mathrm{d}u} + \lambda\dfrac{\sigma^2}{2}} \quad \text{for} \quad \mu \le s \le T. \tag{52}$$

This $\bar{\varphi}_\mu$ is defined in $[\mu, T]$ since $\lambda^{-1} > \dfrac{\sigma^2}{2}$. With this $\bar{\varphi}_\mu$ in mind, we obtain

$$\bar{\psi}_\mu(s) = \psi_\mu(T-s) = \int_\mu^s a(T-v)\bar{\varphi}_\mu(v)\,\mathrm{d}v = \int_\mu^s a(T-v)\varphi_\mu(T-v)\,\mathrm{d}v \quad \text{for} \quad \mu \le s \le T. \tag{53}$$

By substitution $s = T - t$ into (52) and (53), we obtain the desired results (13) and (12). The result (14) is obtained by setting $X_0 = 0$.

□

**A.2 Average and variance of models 1 and 2**

We assume $0 \le t \le 1$ in this section. Here, we assume $T = 1$ following the normalization used in **Section 4**. By elementary calculations, the averages of models 1 and 2 are given by

$$\mathbb{E}[X_t] = \frac{a}{1-c}\left\{(1-t)^c - (1-t)\right\}, \text{ where } c \ne 1 \tag{54}$$

and

$$\mathbb{E}[X_t] = a\frac{1-t}{t+\varepsilon}\left\{(1+\varepsilon)\ln\left(\frac{1}{1-t}\right) - t\right\}, \tag{55}$$

respectively. Likewise, the variances of models 1 and 2 are given by

$$\mathbb{V}\left[X_t\right]=\frac{ac\sigma^2}{1-c}\left\{\frac{1}{c}\left(\left(1-t\right)^c-\left(1-t\right)^{2c}\right)-\frac{1}{1-2c}\left(\left(1-t\right)^{2c}-\left(1-t\right)\right)\right\},\ \text{where}\ \ c\neq 1,2 \tag{56}$$

and

$$\mathbb{V}\left[X_t\right]=a\sigma^2\left(1+\varepsilon\right)\left(\frac{1-t}{t+\varepsilon}\right)^2\frac{1}{1-t}\left\{\left(2+\varepsilon-t\right)\ln\left(\frac{1}{1-t}\right)-\left(2+\varepsilon\right)t\right\}, \tag{57}$$

respectively. Then, the CV is obtained as $\frac{\sqrt{\mathbb{V}\left[X_t\right]}}{\mathbb{E}\left[X_t\right]}$. In model 1, the corresponding results for the cases $c=1,2$ are obtained through the classical L'Hopital's rule.

### A.3 Auxiliary data

**Table A1** shows the total fish count of *P. altivelis* at the study site each year, which is not directly used in the main text. **Figure A1** shows the 10 min fish count data in (a) 2023 and (b) 2024 on a nonlogarithmic scale. **Figure 6** better visualizes the variation in the 10 min migrating fish count, which exceeds four orders of magnitude. **Tables A2 and A3** show the time-average error for the cases corresponding to **Tables 5 and 6**, respectively. These figures show that the time-average error is one order smaller than the maximum one.

**Table A1.** Total fish count of *P. altivelis* at study site each year (courtesy from the Japan Water agency, Nagara River Estuary Barrage Office).

| Year | Total fish count |
|---|---|
| 2024 | 1,236,102 |
| 2023 | 852,596 |
| 2022 | 224,397 |
| 2021 | 403,459 |
| 2020 | 812,342 |
| 2019 | 592,439 |
| 2018 | 847,565 |
| 2017 | 1,171,928 |
| 2016 | 702,028 |
| 2015 | 957,706 |
| 2014 | 608,661 |
| 2013 | 993,089 |
| 2012 | 590,157 |
| 2011 | 841,043 |
| 2010 | 471,415 |
| 2009 | 2,174,478 |
| 2008 | 2,695,955 |
| 2007 | 785,887 |
| 2006 | 130,024 |
| 2005 | 70,157 |
| 2004 | 315,018 |
| 2003 | 437,693 |
| 2002 | 230,000 |
| 2001 | 478,186 |
| 2000 | 568,372 |
| 1999 | 956,441 |
| 1998 | 523,682 |
| 1997 | 534,360 |
| 1996 | 476,319 |
| 1995 | 48,202 |

**Table A2.** Time-average error between computed and theoretical averages for the model 2 in 2023.

| | | Sample size | | |
|---|---|---|---|---|
| | | 10,000 | 100,000 | 1,000,000 |
| Time step | 100 | 3.646.E-04 | 1.329.E-04 | 4.189.E-05 |
| | 1,000 | 1.618.E-04 | 8.445.E-05 | 1.508.E-05 |
| | 10,000 | 2.273.E-04 | 4.749.E-05 | 2.048.E-05 |

**Table A3.** Time-average error between computed and theoretical standard deviations for the model 2 in 2023.

| | | Sample size | | |
|---|---|---|---|---|
| | | 10,000 | 100,000 | 1,000,000 |
| Time step | 100 | 1.676.E-03 | 6.126.E-04 | 8.204.E-04 |
| | 1,000 | 6.380.E-04 | 2.601.E-04 | 8.532.E-05 |
| | 10,000 | 1.327.E-03 | 1.761.E-04 | 1.041.E-04 |

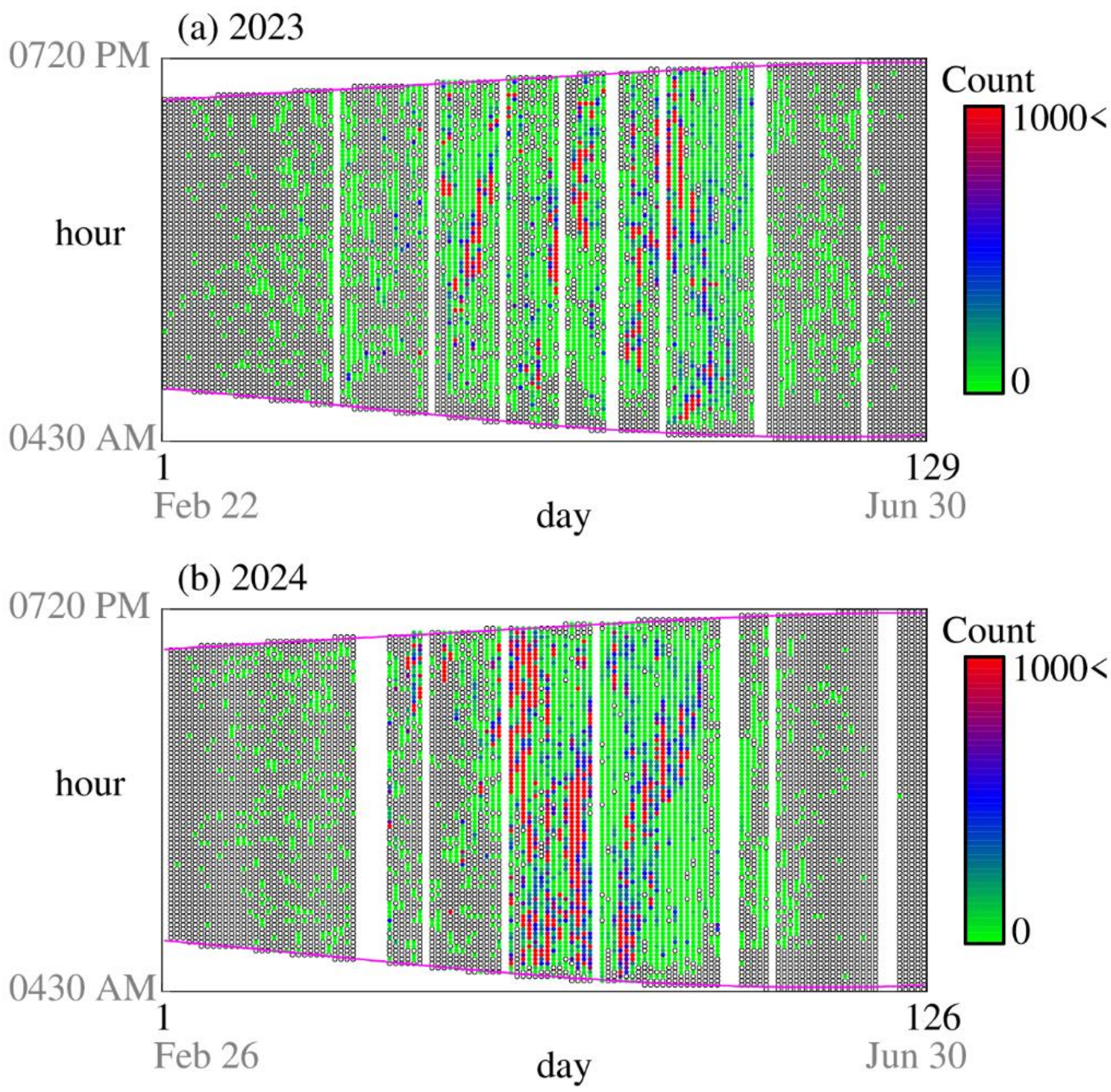

**Figure A1.** 10 min fish count data in (a) 2023 and (b) 2024. Circles are filled if the data are available. Unfilled circles represent a fish count of 0. The magenta curves represent sunrise and sunset times in Nagoya City, as obtained from the National Astronomical Observatory of Japan[6].

[6] https://eco.mtk.nao.ac.jp/koyomi/dni/ (Last accessed on May 15, 2025)